\documentclass[a4paper,11pt, reqno]{amsart}

\usepackage{natbib}
\usepackage{booktabs}
\usepackage{lscape}
\usepackage{float}
\usepackage[table]{xcolor}
\usepackage{array}
\usepackage{booktabs}
\usepackage{amssymb}
\usepackage{amsmath,amsthm}
\usepackage{amsfonts}
\usepackage{enumerate}
\usepackage{graphicx}
\usepackage{hyperref}
\usepackage{color}
\usepackage{cancel}
\usepackage{adjustbox}

\usepackage{geometry}
\usepackage{tikz}
\usetikzlibrary{positioning,fit}

\let\origmaketitle\maketitle
\def\maketitle{
	\begingroup
	\def\uppercasenonmath##1{} 
	\let\MakeUppercase\relax 
	\origmaketitle
	\endgroup
}

\begin{document}

\title[SWIFT]{\Large An Integrated Location-Routing Framework\\for Multi-Type Urban Waste Collection and Recycling Systems}

\author[Blanco, Camacho-Vallejo, \MakeLowercase{and}  Hinojosa]{
{\large V\'ictor Blanco$^{\dagger}$, J. Fernando Camacho-Vallejo$^{\ddagger}$, and Yolanda Hinojosa$^{\star}$}\medskip\\
$^\dagger$Institute of Mathematics (IMAG), Universidad de Granada, Spain\\
$^\ddagger$Tecnologico de Monterrey, Escuela de Ingenieria y Ciencias, Mexico\\
$^\star$Dpt. Applied Economics I, Universidad de Sevilla, Spain\\
\texttt{vblanco@ugr.es}, \texttt{fernando.camacho@tec.mx}, \texttt{yhinojos@us.es}
}

\maketitle

\begin{abstract}
Urban waste collection and recycling systems face increasing operational and environmental challenges due to population growth, heterogeneous waste streams, traffic congestion, and the need for efficient resource recovery. This paper introduces the \textbf{S}ustainable \textbf{W}aste \textbf{I}ntegrated \textbf{F}acility and \textbf{T}ransportation (SWIFT) approach for the design of multi-type urban waste collection and recycling systems. In the proposed system, waste generated at distributed collection points is collected by waste-type-specific vehicles and transported to intermediate consolidation facilities, where it is aggregated before being transferred to treatment plants using larger vehicles. The problem consists of jointly determining the locations of consolidation facilities and treatment plants, together with the associated two-echelon collection and transportation routes for multiple waste streams under a limited investment budget. To address this problem, an integrated location-routing optimization model is developed that simultaneously captures facility-location decisions, waste collection operations, transfer activities, and repeated unloading operations induced by vehicle-capacity limitations. The objective is to minimize the total system cost, including transportation, routing, and handling costs, while satisfying infrastructure investment constraints. Computational experiments based on realistic urban scenarios derived from the city of Medellín, Colombia, demonstrate the benefits of coordinated infrastructure and transportation planning. The results show that strategically located consolidation facilities can improve collection efficiency, enhance vehicle utilization, reduce transportation effort, and support more sustainable urban recycling operations.
\end{abstract}

\section{Introduction}
\label{sec1}

The design of efficient waste collection and recycling systems has become a critical challenge for cities, driven by increasing environmental pressures, regulatory requirements, and economic constraints \citep{nunhes2016evolution, sarra2020optimal}. As part of broader urban sustainability agendas, local governments and private stakeholders are intensifying efforts to promote sustainable waste management practices aligned with circular economy and climate mitigation policies \citep{marrucci2022circular}. Urban areas generate large and heterogeneous volumes of waste as a result of population growth, industrial activity, and changes in consumption patterns. Managing these diverse waste streams poses significant logistical and planning challenges, as each waste type may require specific handling, transportation, and treatment processes \citep{murray2022efficiency}. In this context, the development of advanced decision-support tools for the design and operation of efficient urban waste collection and recycling systems is of paramount importance, particularly for supporting evidence-based policymaking and sustainable urban planning \citep{singh2014progress}.

Within urban waste management systems, routing and transportation decisions play a central role in determining overall operational performance and environmental impact \citep{sahoo2005routing}. In particular, the strategic location of intermediate facilities for waste consolidation can substantially improve the efficiency of the collection process, as highlighted by \cite{ammon2017evaluating}. These facilities allow collection vehicles to offload waste more frequently, thereby increasing routing flexibility and reducing travel distances under vehicle capacity constraints. By jointly considering facility location and routing decisions, urban waste management systems can achieve improvements not only in operational efficiency and cost reduction, but also in energy consumption, fuel usage, and pollutant emissions. Consequently, integrated location-routing approaches contribute directly to the environmental and economic dimensions of urban sustainability.

Motivated by these challenges, this study addresses the design of an integrated urban waste collection network involving multiple waste types through the proposed \textbf{S}ustainable \textbf{W}aste \textbf{I}ntegrated \textbf{F}acility and \textbf{T}ransportation (SWIFT) framework. In this framework, we consider a set of collection points distributed across a city, where different types of waste are generated. Each waste type is initially collected by specialized vehicles that visit all collection points exactly once. Collected waste can be unloaded at strategically located consolidation facilities, which act as intermediate transfer points prior to transportation to treatment plants. For each waste type, a larger vehicle is responsible for transporting waste from consolidation facilities to the corresponding treatment plant. The problem consists of jointly determining the locations of consolidation facilities and treatment plants for each waste type, as well as the associated collection and transportation routes, under a limited investment budget. The objective is to minimize the total system cost, including transportation, routing, and handling costs, while facility location costs are constrained by the available investment budget.

\begin{figure}
\begin{center}
\includegraphics[width=1\textwidth]{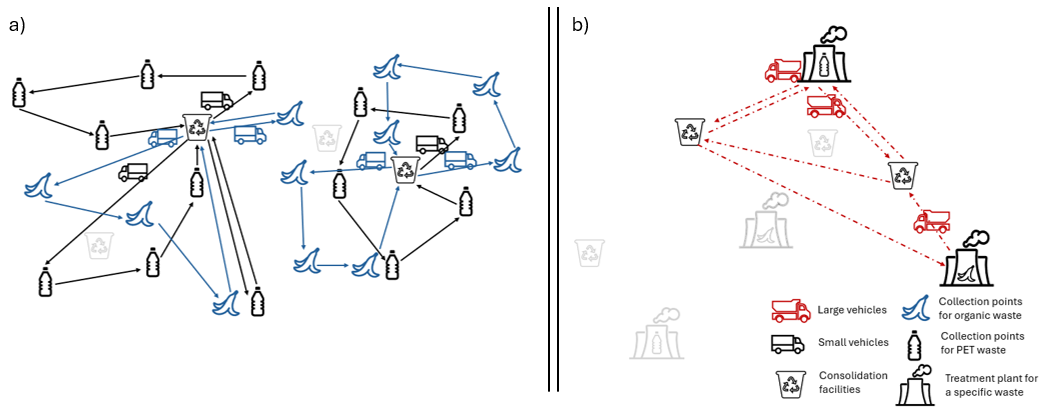}
\end{center}
\caption{Illustration of the two operational phases of the SWIFT framework. \label{fig:sketch}}
\end{figure}

To illustrate the problem under study, consider the schematic representation shown in Figure \ref{fig:sketch}a), which depicts the first phase of the system and considers two types of urban waste: PET bottles and organic waste. In this phase, the locations of consolidation facilities are determined to support the routing operations performed by small collection vehicles. For each waste type, all collection points must be visited exactly once, and the use of strategically located consolidation facilities allows PET bottles and organic waste to be efficiently aggregated. Small vehicles are not required to start their routes from a central depot; instead, routes originate at consolidation facilities, which may be visited multiple times due to vehicle capacity limitations. Additionally, multiple routes may be operated as needed to ensure full coverage of the urban area.

The second phase of the system is illustrated in Figure \ref{fig:sketch}b). In this phase, treatment plants specialized in processing each waste type are located, and once the waste has been consolidated, larger vehicles are dispatched to collect the aggregated waste from the selected consolidation facilities and transport it to the corresponding treatment plants.

\begin{figure}
\begin{center}
\includegraphics[width=0.6\textwidth]{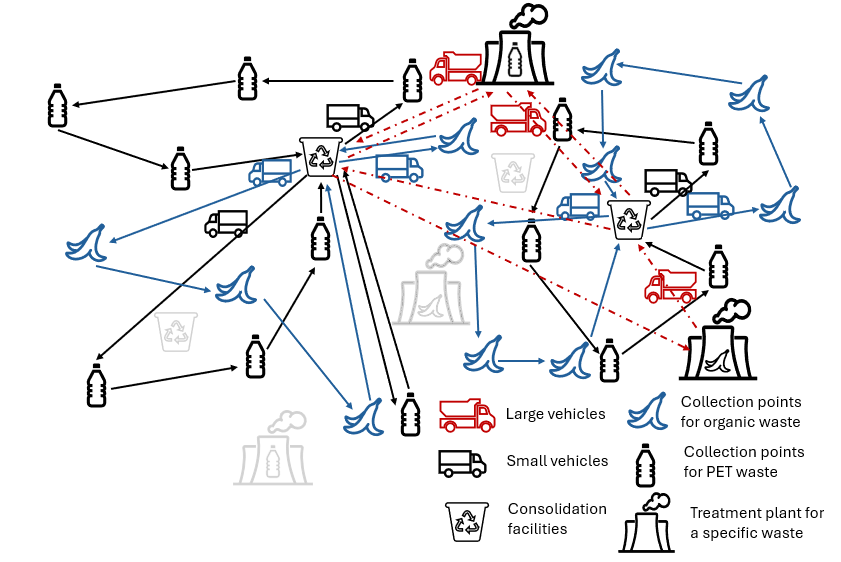}
\end{center}
\caption{Example of the SWIFT framework. \label{fig:integrating}}
\end{figure}

This two-phase structure reflects typical operational practices in urban waste management systems and highlights the need for coordinated location and routing decisions across different system levels. It is worth noting that, although the system can be conceptually decomposed into two phases for illustrative purposes, all decisions are jointly optimized within a single integrated framework, resulting in a scheme similar to that presented in Figure \ref{fig:integrating}.

One of the advantages of the proposed model is its versatility, as it can be applied to a wide range of practical contexts. For instance, in municipal or urban recycling programs, consolidation facilities may operate as transfer or sorting stations, enabling more flexible routing and reducing the number of long-distance trips required for waste transportation (see \cite{asefi2015mathematical, liao2025learning}). As noted earlier, private waste management stakeholders may also benefit from improved facility utilization and reduced operational costs by integrating consolidation points into their service networks. These application contexts contribute to more sustainable urban logistics by optimizing the use of available resources, minimizing environmental impacts, and enhancing the efficiency of recycling operations, as discussed by \cite{purkayastha2015collection}.

The remainder of this paper is organized as follows. Section~\ref{sec:2} reviews the related literature and positions the proposed SWIFT framework within existing waste collection and location-routing studies. Section~\ref{sec:3} introduces the SWIFT problem and presents the proposed mathematical formulation. Additionally, a procedure for obtaining an initial feasible solution is presented to support the solution process. Section~\ref{sec:4} reports the computational experiments based on urban scenarios derived from the city of Medell\'in and discusses the obtained managerial and sustainability insights. Finally, Section~\ref{sec:5} summarizes the main conclusions and outlines directions for future research.

\section{Literature Review}\label{sec:2}

The literature review is organized into three main related topics. First, collection routing problems in the context of waste management are discussed. Next, the use of intermediate facilities within waste collection systems is summarized. Finally, two-echelon location-routing problems are detailed to contextualize our problem within the existing literature. 

\subsection{Waste Collection Location-Routing Problems}

Waste collection and recycling systems have received considerable attention in the literature due to their economic, environmental, and social relevance \citep{pires2019sustainable}. Early studies focused on routing decisions for municipal solid waste collection, primarily addressing operational aspects such as vehicle capacities \citep{akhtar2017backtracking, hannan2018capacitated}, vehicle characteristics \citep{zhang2026optimization}, service frequencies \citep{angelelli2002application, hemmelmayr2014models}, route duration limits \citep{benjamin2010metaheuristics}, and time windows \citep{tung2000vehicle, kim2006waste}. Over time, these models evolved to incorporate additional characteristics, including heterogeneous fleets, stochastic waste generation, dynamic routing conditions, and sustainability-oriented objectives, as discussed in the survey by \cite{belien2014municipal} and more recently highlighted by \cite{hess2024waste} and \cite{li2025literature}.

In parallel, another stream of research has investigated strategic decisions related to the location of waste management infrastructure, including recycling centers, transfer stations, treatment facilities, and landfills; see \cite{adeleke2020facility} for a review. These studies typically aim to determine cost-effective facility configurations, as in \cite{asefi2015mathematical}, while satisfying capacity and environmental constraints (e.g., \citep{farrokhi2017metaheuristics, rabbani2017solving}). More recently, the strong interdependence between facility location and routing decisions has motivated the development of integrated location-routing models for waste management systems. \cite{maheo2023solving} discusses the advantages of solving integrated systems rather than adopting the traditional sequential approach. As noted several years ago by \cite{min1998combined}, integrated approaches can lead to more efficient networks, which nowadays may also contribute to the development of more sustainable recycling systems.

Within this context, researchers have addressed recycling systems involving different types of waste materials, such as plastics, glass, paper, electronic waste, medical, and organic waste \citep{BHZ2024,mofid2019multi, tang2025improved, rodriguez2026bi}. \cite{rodrigues2016waste} noted that multi-type systems introduce additional complexity due to differentiated handling requirements, specialized vehicle fleets, and waste-specific treatment processes. Although the literature on waste collection location-routing problems has grown substantially (see the previously mentioned recent reviews), relatively limited attention has been devoted to integrated urban recycling systems that simultaneously consider multiple waste types, routing operations, and the location of both consolidation facilities and treatment plants. 

\subsection{Considering Intermediate Facilities for Waste Collection}

Intermediate facilities have played an important role in the design of waste collection and recycling systems due to their ability to improve transportation efficiency and reduce operational costs \citep{asefi2015mathematical}. As noted by \cite{wei2024multi}, in practice, these facilities may operate as transfer stations, consolidation centers, or temporary storage locations (huts) where collected waste is aggregated before being transported to final treatment or disposal facilities. The use of intermediate facilities may substantially improve existing waste collection plans by enabling more flexible routing structures and reducing transportation inefficiencies. As observed in practical applications such as the study by \cite{taverna2025enhancing}, relatively small operational modifications can generate important cost savings while simultaneously improving the robustness and practicality of collection operations. Recently, \cite{albareda2025coordinating} further emphasized the relevance of jointly coordinating intermediate drop-off locations and collection routes under budget limitations, highlighting the operational and economic benefits that integrated routing and consolidation decisions may provide in transportation and collection systems.

The incorporation of intermediate facilities also introduces additional operational complexity into the routing process, since collection vehicles may perform multiple unloading operations during a route and waste flows must be coordinated with downstream transportation activities \citep{lu2015smart}. \cite{markov2016integrating} studied a vehicle routing problem with intermediate facilities incorporating heterogeneous fleets and flexible assignment of destination depots, motivated by waste collection operations involving long tours and repeated unloading activities; however, facility location decisions are not considered. \cite{asefi2015mathematical} proposed an integrated solid waste management model that considers transfer stations, recycling facilities, treatment centers, and disposal sites within a multi-echelon transportation network. Their framework incorporates multiple waste types, heterogeneous vehicle fleets, compatibility constraints, and time windows under both deterministic and stochastic settings. However, explicit collection routes visiting waste generation points are not considered. Instead, transportation is modeled as direct trips between waste generation points and the different facilities within the network. 

 Despite these advances, relatively limited attention has been devoted to integrated urban recycling systems that simultaneously consider multiple waste types, specialized vehicle fleets, explicit collection routes through waste generation points, repeated unloading operations at intermediate consolidation facilities, and the location of both consolidation facilities and treatment plants. The proposed SWIFT framework contributes to this literature by integrating these decisions within a unified two-echelon optimization model, in which waste-type-specific collection routes, intermediate unloading operations, facility-location decisions, and downstream transportation activities are optimized in a coordinated manner.

\subsection{An Overview on Two-Echelon Location-Routing Problems}

Two-echelon location-routing problems (2E-LRPs) represent hierarchical transportation systems commonly arising in urban logistics and distribution networks. In these problems, goods or commodities are transported through two interconnected transportation levels: a first echelon involving long-haul transportation between central facilities and intermediate satellites, and a second echelon responsible for local distribution or collection operations \citep{cuda2015survey, darvish2019flexible}. According to \cite{cheng2022multi}, this structure allows waste collection and transfer activities to be coordinated across different vehicle types and service regions, improving operational efficiency and reducing transportation effort.

In the context of waste collection systems, \cite{yu2020two} proposed a generalized 2E-LRP with multiple objectives that simultaneously considers intermediate collection centers and final disposal facilities. Their approach incorporates routing, facility, and vehicle costs, allowing the model to represent different operational and sustainability objectives arising in waste collection applications. Nevertheless, their formulation assumes a restrictive routing structure in which each route starts and ends at the same facility, customers are assigned to specific facilities, multiple vehicles may be required to serve the assigned customers, and only one waste type is considered. In a broader integrated waste-management context, \cite{asefi2019integrated} proposed a tri-echelon logistics network that simultaneously considers transfer stations, recycling facilities, treatment centers, and disposal sites under uncertainty. Their framework incorporates multiple waste types, heterogeneous fleets, and detailed routing decisions within a complete waste-management chain. However, facility location decisions are not considered.

In contrast, the proposed SWIFT framework jointly optimizes facility location, treatment plant selection, waste-type-specific collection routes, repeated unloading operations, and downstream transportation flows within a two-echelon urban recycling structure. Each waste type is served by a single specialized vehicle whose route may start and end at different consolidation facilities and unload as many times as necessary. Thus, SWIFT integrates multiple waste streams, location decisions, and coordinated collection and transfer operations within a unified optimization framework.

\subsection{Research Gap and Contributions}

Recent urban sustainability policies increasingly emphasize the need for integrated recycling systems capable of reducing operational costs, transportation effort, and environmental impacts. However, most existing waste collection studies still treat facility location, routing, and transfer operations as independent planning problems. Moreover, the recent review by \citet{hess2024waste} highlights the integration of pick-up and delivery interactions within waste collection systems as an important and still underexplored research direction.

From a practical perspective, such integration is essential to avoid suboptimal decisions arising from the separate optimization of interconnected system components. Decisions regarding facility placement directly affect routing patterns and transfer operations, while routing choices influence facility utilization and capacity requirements. Ignoring these interactions may lead to inefficient waste collection systems characterized by higher operational costs, excessive travel distances, increased environmental burdens, and poor service quality. In turn, these shortcomings can reduce the acceptance and satisfaction of both citizens and municipal organizations, compromising the long-term viability and utilization of the infrastructure. Consequently, integrated planning approaches are necessary to ensure that investments in recycling and waste management systems effectively contribute to sustainability, environmental protection, and resource efficiency objectives.

Motivated by these challenges and by the need to capture the strong interdependencies among facility location, waste collection, transfer operations, and multiple waste streams, this paper introduces the SWIFT framework, an integrated optimization approach for urban recycling systems with intermediate consolidation facilities. By jointly optimizing all relevant strategic and operational decisions, the proposed framework aims to improve system efficiency, reduce operational (and then,  environmental) costs, and support the design of sustainable and economically viable recycling networks. The main contributions of this work are summarized as follows:
\begin{itemize}
\item[-] We propose a novel multi-type location--routing model that jointly determines the locations of consolidation facilities and treatment plants, together with the associated collection and transfer routes for different waste streams. By integrating these interdependent decisions, the model avoids the inefficiencies associated with sequential planning approaches and enables the design of more effective recycling networks.
\item[-] The proposed framework integrates collection, unloading, consolidation, and downstream transportation decisions within a unified two-echelon structure, explicitly accounting for vehicle-capacity limitations and repeated unloading operations. This integrated perspective allows for a more efficient utilization of resources and infrastructure throughout the recycling process.
\item[-] The model captures key operational characteristics of urban recycling systems, including specialized vehicle fleets and coordinated interactions between collection and transfer activities. These features enhance the realism of the planning process and support implementable solutions for complex urban environments.
\item[-] From a sustainability perspective, the proposed approach promotes the design of recycling systems that reduce unnecessary transportation effort, improve vehicle utilization, and lower fuel consumption. These improvements contribute not only to lower operational costs but also to reduced greenhouse-gas emissions, congestion, and environmental impacts associated with waste collection activities.
\item[-] Computational experiments based on realistic scenarios derived from the city of Medell\'in provide managerial insights into the value of intermediate consolidation facilities. The results show how integrated planning can improve routing efficiency, reduce transportation requirements, and support more sustainable and economically viable urban recycling operations.
\end{itemize}

\section{Mathematical Formulation}\label{sec:3}

In this study, we address the problem of designing an integrated waste
collection network for recycling purposes. We consider a set of collection
points at which multiple types of waste are generated, with each waste type
collected by specialized vehicles. To enhance the operational efficiency
of the collection process, consolidation facilities can be strategically located
along the routes. The set of located consolidation facilities is shared by all
waste types, whereas collection decisions are defined separately for each waste
type. These facilities enable vehicles to unload accumulated waste, thereby
reducing routing distances and overall transportation costs. 
For each waste type, a single small specialized vehicle is assumed to visit all
associated collection points exactly once, following a route that may start and
end at different consolidation facilities. Nevertheless, due to vehicle
capacity limitations, each consolidation facility may be visited multiple times
to offload waste collected from the corresponding collection points.

Once the waste has been collected from the collection points and unloaded at the consolidation facilities by the corresponding small vehicles, it must be transported to the treatment plants. Unlike consolidation facilities, treatment plants are waste-type specific: for each waste type, a single treatment plant is established, and these plants may be located at different sites. A large vehicle is assigned to collect the accumulated waste from the selected consolidation facilities and transport it to the corresponding treatment plant. The route of each large vehicle starts and ends at its designated treatment plant. However, due to vehicle-capacity limitations, the vehicle may need to return to the treatment plant multiple times during its route to unload the waste currently being carried. Additionally, since each consolidation facility must be visited exactly once by the corresponding large vehicle during the second phase, this constraint ensures that the maximum storage capacity of a consolidation facility for a given waste type does not exceed the capacity of that vehicle.

Our objective is, given a fixed budget for establishing both consolidation facilities and treatment plants for each waste type, to determine their optimal locations and design the corresponding transportation plan that transfers all waste from the collection points to the appropriate treatment plants via the consolidation facilities. The goal is to minimize the total system cost, which includes transportation, routing, and handling costs.

To address this problem, we propose a mixed-integer linear programming (MILP) model that captures the key decisions involved in the design and operation of the waste collection network. The model jointly determines the locations of consolidation facilities and treatment plants, the routing of small and large vehicles, and the flow of waste across the different network levels, while explicitly accounting for capacity constraints and budget limitations.

The index sets and parameters used in our model are summarized in Table~\ref{tab:params}.

\begin{table}[H]
  \centering
  \adjustbox{scale=0.88}{
  \begin{tabular}{@{} >{\ttfamily}l @{\hspace{1em}} l @{}}
    \multicolumn{2}{@{}l}{\textbf{\textsf{Sets}}} \\
    \cmidrule[0.8pt](l{5pt}r{5pt}){1-2}
    $I$       & Collection points. \\
    $J$       & Potential sites for locating consolidation facilities. \\
    $L$ & Types of waste. \\ 
    $K_{\ell}$ & Potential sites for locating the treatment plant for waste type $\ell \in L$.\\
    \multicolumn{2}{@{}l}{\textbf{\textsf{Parameters}}} \\
    \cmidrule[0.8pt](l{5pt}r{5pt}){1-2}
    $d_{i\ell}$     & Amount of waste of type $\ell \in L$ to be collected from the collection point $i \in I$. \\
   $c_{ij\ell}^1$ & Operational costs associated with the small vehicle transporting waste of type\\
    & $\ell \in L$ from  $i\in I\cup J$ to  $j\in I\cup J$.\\
     $c_{jk\ell}^2$ &  Operational costs associated with the large vehicle transporting waste of type\\
    & $\ell \in L$ from  $j\in J\cup K_\ell$ to  $k\in J\cup K_\ell$.\\
    $f^1_j$ & Fixed cost for locating consolidation facility $j\in J$. \\
    $f^2_{k\ell}$ & Fixed cost for locating a treatment plant $k\in K_{\ell}$ for type of waste $\ell \in L$. \\
    $B$ &  Available budget for locating consolidation facilities and treatment plants.\\
    $\rho^1_{\ell}$ & Capacity of the small vehicle collecting waste of type $\ell \in L$.\\
   $\rho^2_{\ell}$ & Capacity of the large vehicle collecting waste of type $\ell \in L$.\\
 \hline
  \end{tabular}}
  \caption{Index sets and parameters of the SWIFT model.}
  \label{tab:params}
\end{table}

The variables of our model represent the main decisions made in the SWIFT framework, namely facility location, routing, and waste flow 
decisions across the recycling network. These decision variables are summarized in Table~\ref{tab:vars}. 

\begin{table}[tbh]
  \centering
  \adjustbox{scale=0.85}{
  \begin{tabular}{@{} >{\ttfamily}l @{\hspace{1em}} l @{}}
    \multicolumn{2}{@{}l}{\textbf{\textsf{ Decision variables.}}} \\
    \cmidrule[0.8pt](l{5pt}r{5pt}){1-2}
    $y_{j}$ & Binary variable that is 1 if consolidation facility $j \in J$ is located;\\
    & and 0 otherwise.\\
    $\xi_{k\ell}$ &   Binary variable that is 1 if treatment plant $k \in K_{\ell}$ is located for type \\
    & of waste $\ell \in L$; and 0 otherwise.\\
    $x_{ij_\ell}$ & 
      Binary variable equal to 1 if the small vehicle collecting waste of type  $\ell \in L$ \\
      & traverses the arc from $i\in I\cup J$ to $j\in I\cup J$ ($i\neq j$); and 0 otherwise.\\
    $z_{jk\ell}$ & Binary variable equal to 1 if the large vehicle collecting waste of type  $\ell \in L$ \\
      & traverses the arc from  $j\in J\cup K_{\ell}$ to $k\in J\cup K_{\ell}$ ($j\neq k$); and 0 otherwise.\\
    \multicolumn{2}{@{}l}{\textbf{\textsf{ Auxiliary variables}}} \\
    \cmidrule[0.8pt](l{5pt}r{5pt}){1-2}
    $u_{i\ell}$ & Quantity of waste of type $\ell \in L$  transported by the small vehicle upon departing \\
    &from $i \in I \cup J$.\\ 
    $h_{ij\ell}$ & Quantity of waste of type $\ell \in L$ unloaded in $j\in J$ arriving from $i \in I$. \\
    $w_{j\ell}$ &  Quantity of waste of type $\ell \in L$ transported by the large vehicle upon departing\\
    & from $j \in J\cup K_\ell$.\\\hline
  \end{tabular}}
  \caption{Variables used in the SWIFT model.}
  \label{tab:vars}
\end{table}

Routing variables $x_{ij\ell}$ are defined so as to reflect the practical operation of the waste collection system. In particular, once a small vehicle reaches a consolidation facility, it unloads all the waste it has collected for the corresponding waste type. As a result, there is no operational reason for a small vehicle to travel directly from one consolidation facility to another without first visiting a collection point, and such movements are therefore excluded from the model. This is imposed by $x_{ij\ell} = 0$ for $i \in J$ and $j \in J$.

Similarly, for each waste type, only one treatment plant is selected being the final destination for all collected waste. Then,  $z_{jk\ell} = 0$ for $j \in K_\ell$ and $k \in K_\ell$.

In what follows, we present the objective function and describe the constraints that define the proposed optimization model.

\subsection{Objective Function}
The objective of the proposed MILP model is to minimize the total system cost, which includes all transportation 
and handling costs associated with the collection and transfer of different waste types. To achieve this, the objective function is expressed as follows: 

\begin{equation}
\min \quad
\sum_{\ell \in L}\sum_{i,j\in I\cup J}c^1_{ij\ell}, x_{ij\ell}
+
\sum_{\ell \in L}\sum_{j,k\in J\cup K_\ell}c^2_{jk\ell}, z_{jk\ell}
\end{equation}

The first term represents the operational costs associated with the small vehicles, which transport waste of type $\ell \in L$ along arcs connecting nodes in $I \cup J$. The second term captures the operational costs of the large vehicles, which transport consolidated waste of type $\ell \in L$ between consolidation facilities and treatment plants in $J \cup K_\ell$.

The operational costs $c_{ij\ell}^1$ and $c_{jk\ell}^2$ may include both transportation 
and handling components. 
When only transportation costs are considered, they are directly associated with the routing decisions represented by the binary variables $x_{ij\ell}$ and $z_{jk\ell}$. If handling costs related to waste unloading are also included, the model must additionally account for visits to consolidation facilities and treatment plants. These visits can likewise be inferred from the routing variables $x_{ij\ell}$ and $z_{jk\ell}$.

In summary, the objective function integrates all routing decisions for both small and large vehicles and quantifies their associated costs. By minimizing this function, the model identifies the set of transportation routes and waste flows that result in the lowest total operational cost for the system.

\subsection{Constraints}

The constraints presented below define the feasible region of the proposed MILP model. They formalize the logical, operational, and structural relationships among the previously introduced parameters and decision variables. Together, these constraints ensure that facility-location decisions, routing choices, and waste-flow assignments are internally consistent and satisfy the technical requirements of the system, including budget limitations, facility locating decisions, vehicle routing rules, and minimum service levels at consolidation facilities.\\

\noindent\textbf{$\star$  Budget constraint}

The total investment required to locate consolidation facilities and treatment plants must not exceed the available budget:
\begin{equation}
    \sum_{j \in J} f^1_j y_j \;+\; \sum_{\ell\in L} \sum_{k\in K_\ell} f^2_{k\ell} \xi_{k\ell} \leq B.
    \label{budget}
\end{equation}
This constraint aggregates all fixed installation costs. The first term (sum of $f^1_j y_j$) represents the cost of installing the consolidation facilities, whereas the second term (sum of $f^2_{k\ell}\xi_{k\ell}$) accounts for the cost of locating the treatment plants for each waste type. The inequality enforces that the total investment does not exceed the available budget $B$.

Note that in the case where all consolidation facilities have the same fixed set-up cost, that is, 
$f^1_j = F^1$ for all $j \in J$, and where, for each waste type $\ell \in L$, all treatment 
plants incur an identical set-up cost, that is, $f^2_{k\ell} = F^2_{\ell}$ for all 
$k \in K_\ell$, the budget constraint simplifies considerably: 
\begin{align}
    \sum_{j \in J} y_j 
    \leq 
    p,  \quad \text{with} \quad  p= 
    \left\lceil 
        \frac{B - \sum_{\ell \in L} F^2_{\ell}}
        {F^1}
    \right\rceil,
    \label{mC:ctr2b}
\end{align}
providing an upper bound on the number of consolidation facilities that may be opened 
within the available budget.

Another useful observation is that a lower bound on the number of consolidation facilities that must be open to accommodate all waste types is given by:
\begin{equation}
\sum_{j \in J} y_j \geq N, \quad \text{with} \quad N = \max_\ell \left\lceil \frac{\sum_{i \in I} d_{i\ell}}{\rho^2_\ell} \right\rceil.
\end{equation}
The bound arises from the fact that the largest vehicle used in the second phase has a limited capacity $\rho^2_\ell$ and can visit each consolidation facility only once. Consequently, at least $N$ facilities must be open to ensure that all collected waste of each type can be transported without exceeding the vehicle capacity. \\

\noindent\emph{- First Transportation Phase:}\\

\noindent\textbf{$\star$  Activation of consolidation facilities}

Routing arcs involving a consolidation facility $j \in J$ can only be used if the facility is located:
\begin{align}
        x_{ij\ell} &\leq y_j \;\; \text{and} \ \  x_{ji\ell} \leq y_j, & \forall i \in I,\; j \in J,\; \ell \in L,\label{const_act1}
\end{align}
These inequalities enforce that no incoming or outgoing arc associated with facility $j \in J$ can be used unless $y_j = 1$.\\

\noindent\textbf{$\star$  Routing at collection points}

For each collection point and each waste type, there must be exactly one outgoing and one incoming arc:
\begin{align}
     \sum_{j \in I\cup J} x_{ij\ell} &=\sum_{j \in I\cup J} x_{ji\ell} = 1, \; & \forall i \in I,\; \ell \in L, \label{const_routing}
\end{align}
These constraints ensure that each collection point is visited exactly once by a small vehicle for each waste type, under the assumption that small vehicle capacity exceeds the demand at collection points.\\

\noindent\textbf{$\star$  Minimum service at consolidation facilities}
A consolidation facility is considered operational for a waste type only if it is
actually visited by the corresponding collection route. To enforce this
condition, the following constraints are imposed:
\begin{align}
\sum_{i\in I} x_{ij\ell} \geq y_j, 
\quad &\forall j\in J,\ \ell\in L \label{const_minservice}
\end{align}
These constraints ensure that, whenever a consolidation facility $j$ is located
(i.e., $y_j=1$), it receives at least one incoming routing arc for each waste
type $\ell\in L$, preventing the location of
facilities that are not visited by any collection route,  even if sufficient budget is available.\\

\noindent\textbf{$\star$  Flow conservation constraints}

The following constraints ensure the correct propagation of waste flows along the routes followed by the small vehicles.
\begin{align}
     u_{i\ell} + d_{j\ell} - 2\rho^1_{\ell}\big(1 - x_{ij\ell}\big) \leq & u_{j\ell} \leq u_{i\ell} + d_{j\ell} + \rho^1_{\ell}\big(1 - x_{ij\ell}\big),
     && \forall i \in I \cup J,\; j \in I,\; \ell \in L,
     \label{ctr:flow1} \\
     u_{i\ell} - \rho^1_{\ell}\big(1 - x_{ij\ell}\big) \leq & h_{ij\ell} \leq u_{i\ell} + \rho^1_{\ell}\big(1 - x_{ij\ell}\big),
     && \forall i \in I,\; j \in J,\; \ell \in L,
     \label{ctr:flow4} \\
     h_{ij\ell} &\leq \rho^1_{\ell}\, x_{ij\ell},
     && \forall i \in I,\; j \in J,\; \ell \in L.
     \label{ctr:flow6}\\
          u_{j\ell}& = 0, 
     && \forall j \in J,\; \ell \in L.
     \label{ctr:flow3}\\
     \sum_{i\in I} h_{ij\ell} &\leq \rho^2_\ell y_j, 
     && \forall j\in J, \; \ell \in L. \label{const_capacity}
\end{align}
Constraints \eqref{ctr:flow1} impose load continuity along the arcs visited by the small vehicle. When $(i,j)$ is traversed for waste type $\ell$ ($x_{ij\ell}=1$), these constraints reduce to $
u_{j\ell} = u_{i\ell} + d_{j\ell},
$
ensuring that the vehicle load upon leaving $j$ equals the load arriving from $i$ plus the amount collected at $j$. When the arc is not used, the terms involving $\rho^1_{\ell}$ relax the constraints. 
Constraints \eqref{ctr:flow4} guarantee consistency between the load carried upon leaving node $i$ and the amount unloaded at consolidation facility $j$. If arc $(i,j)$ is used then the constraints imply
$
h_{ij\ell} = u_{i\ell},
$
meaning that the entire load carried is unloaded at the consolidation facility. Constraints \eqref{ctr:flow6} enforce that unloading can occur only on arcs that are actually traversed by the vehicle and that the unloaded amount does not exceed the vehicle capacity. Constraints \eqref{ctr:flow3} force the load to be zero when leaving any consolidation facility because all waste is unloaded upon arrival at such facilities. Finally, constraints \eqref{const_capacity} directly link the admissible inflow to each facility to the second phase vehicle capacity ($\rho^2_\ell$), ensuring that the accumulated waste at a facility does not exceed the maximum quantity  that can be collected in a single visit by the large vehicle during the second phase, while enforcing that no unloading occurs at facilities that are not open ($y_j = 0$).
 Together, these constraints enforce consistency among routing decisions, collected loads, and unloading operations.\\




\noindent\textbf{$\star$  Route continuity constraints at consolidation facilities}

For each waste type, the collection process is carried out by a single vehicle that must serve all collection points and the consolidation facilities selected in the network. The routing decisions are therefore required to ensure route continuity, meaning that the vehicle follows a single, uninterrupted path through the network.
Operationally, this is enforced by requiring that every time the vehicle enters a node, it also leaves that node, except possibly at the initial and final locations of the route.

In the proposed formulation, constraints \eqref{const_routing} ensure that, at every collection point, the number of times the vehicle arrives is equal to the number of times it departs (one time), thereby guaranteeing a consistent visit of all collection locations. At consolidation facilities, this balance is enforced only for those facilities that are selected to operate. Moreover, among the located consolidation facilities, at most two are allowed to present an imbalance between arrivals and departures, corresponding to the initial and final locations of the collection route.

To model these admissible imbalances 
we introduce binary variables that indicate whether a consolidation facility serves as the starting or ending point of the route for a given waste type:\\

{\small\begin{tabular}{@{} >{\ttfamily}l @{\hspace{1em}} l @{}}
    \cmidrule[0.8pt](l{5pt}r{5pt}){1-2}
    $e_{j\ell}$ & 
      Binary variable equal to 1 if if facility $j \in J$ receives the vehicle one  \\
      &  more time than it sends it out for waste type $\ell\in L$; and 0 otherwise.\\
      &\\
    $s_{j\ell}$ &  
    Binary variable equal to 1 if if facility $j\in J$  sends out the vehicle one  \\
      &  more time than it receives it  for waste type $\ell\in L$; and 0 otherwise.\\
      \hline
\end{tabular}}\\

These variables are directly linked to the routing decisions through
\begin{equation}
    \sum_{i\in I} x_{ji\ell} - \sum_{i\in I} x_{ij\ell} = s_{j\ell} - e_{j\ell},
    \qquad \forall j\in J,\; \ell\in L,
    \label{ctr:euler1}
\end{equation}
which ensures that $s_{j\ell}=1$ identifies facility $j \in J$ as the starting point of the route for waste type $\ell \in L$, while $e_{j\ell}=1$ identifies it as the ending point.

To ensure that, for each waste type, a route has at most one starting facility and at most one ending facility, we impose the following constraints:
\begin{align}
    \sum_{j\in J} e_{j\ell} \leq 1 \ \ & \text{and} \ \ \sum_{j\in J} s_{j\ell}\leq 1, & \forall \ell \in L. \label{ctr:euler2b}
\end{align}

Finally, these arrival-departure imbalances are only allowed at consolidation facilities that are actually opened:
\begin{equation}
e_{j\ell} + s_{j\ell} \le y_j, \quad \forall j \in J,\; \ell \in L.
\end{equation}

Additionally, to guarantee that the routes in this first phase are connected, we impose the following standard cutset constraints:
\begin{equation}
\sum_{i\in S,\; j\in S^c} x_{ij\ell} 
\;\geq\; y_r + y_{r'} - 1 - \sum_{k\in S^c} s_{k\ell},
\quad 
\forall S\subset I\cup J,\; r\in S\cap J,\; r'\in S^c\cap J,\; \ell\in L.
\label{ctr:conn1}
\end{equation}

Constraints \eqref{ctr:conn1} ensure that, for any subset of nodes $S\subset I\cup J$ and its complement $S^c$, there is at least one active arc connecting a node in $S$ with a node in $S^c$.  

These constraints are only active when $\sum_{k\in S^c} s_{k\ell}=0$, $y_r=1$, and $y_{r'}=1$. Intuitively, they examine all possible partitions of the node set and enforce connectivity across each cut $(S, S^c)$. In particular, if the starting point of the route is not in $S^c$ ($\sum_{k\in S^c} s_{k\ell}=0$), either because it lies in $S$ or because the route forms a cycle without a designated starting point, then there must exist at least one arc leaving $S$ and entering $S^c$.  

To limit the number of constraints, only subsets $S$ for which both $S$ and its complement $S^c$ contain at least one located consolidation facility (i.e., $y_r = 1$ and $y_{r'} = 1$) are considered. This restriction is justified because the flow conservation constraints already ensure that all collection points are connected to at least one located consolidation facility for waste unloading. Consequently, subsets composed exclusively of collection points need not be considered, as they cannot be disconnected from facilities capable of receiving waste without violating flow feasibility.  

In this way, constraints \eqref{ctr:conn1} prevent the formation of disconnected routes and ensure that no subset of nodes containing located consolidation facilities becomes isolated from the rest of the network. \\

\noindent \emph{- Second transportation phase:}\\

\noindent\textbf{$\star$  Treatment plant selection}

For each waste type $\ell$, exactly one location from the candidate set $K_\ell$ is selected to establish the corresponding treatment plant. This is enforced by the following constraints:
\begin{equation}
    \sum_{k \in K_\ell} \xi_{k\ell} = 1, \qquad \forall \ell \in L.
\end{equation}

\noindent\textbf{$\star$  Restrictions on connections to located treatment plants}

 Arcs incident to non-located treatment plants are forbidden. This is captured through:
 \begin{align}
    z_{kj\ell} &\leq \xi_{k\ell} 
\;\; \text{and} \ \  z_{jk\ell} \leq \xi_{k\ell},
    & \forall j\in J,\; \ell\in L, \; k\in K_\ell.
    \label{ctr:noClosed2}
\end{align}
These constraints prevent any incoming or outgoing flow involving a treatment plant $k \in K_\ell$ that is not open $( \xi_{k\ell}=0)$. \\



\noindent\textbf{$\star$  Route continuity}

In the second phase, the waste accumulated at the located consolidation facilities is
transported by a large vehicle to the treatment plant selected for each waste type. In this phase, the vehicle follows a closed route that starts and ends at the same location, ensuring a continuous transportation process. This requires suitable balance conditions at both consolidation facilities and treatment plants.  

For each waste type $\ell \in L$, every located consolidation facility must be visited exactly once by the large vehicle. This is enforced by requiring exactly one incoming and one outgoing arc at each open facility:

\begin{align}
    \sum_{j\in J\cup K_\ell} z_{jk\ell} \;  = \;\sum_{j\in J\cup K_\ell} z_{kj\ell}  &= y_k, && \forall k\in J,\; \ell\in L, \label{ctr:eulerCycleA}
\end{align}
At the treatment plant, route continuity is ensured by requiring the number of incoming
arcs to equal the number of outgoing arcs:
\begin{equation}
    \sum_{j\in J} z_{jk\ell} = \sum_{j\in J} z_{kj\ell}, 
    \qquad \forall \ell\in L, \; k\in K_\ell. \label{ctr:eulerCycleTP}
\end{equation}

Together, these constraints define a continuous closed route for each waste type,  which each located consolidation facility is served exactly once while the treatment plant may be visited multiple times if required. \\

\noindent\textbf{$\star$  Flow conservation} 

Let $   H_{j\ell} = \sum_{i\in I} h_{ij\ell} $ be an auxiliary variable representing the total amount of waste of type $\ell\in L$ delivered and unloaded at consolidation facility $j \in J$. 

To ensure a consistent evolution of the vehicle load along the selected routes, the following constraints are imposed

\begin{align}
    w_{j\ell} + H_{k\ell} - 2\rho^2_\ell (1 - z_{jk\ell})\leq w_{k\ell} &\leq w_{j\ell} + H_{k\ell} + \rho^2_\ell (1 - z_{jk\ell}),
    && \forall \ell\in L,\; k\in J,\; j\in J\cup K_\ell. 
    \label{ctr:flow2phaseB}\\
    w_{k\ell} &= 0, && \forall  \ell\in L, \; k \in K_\ell. \label{ctr:flow2phase0}\\
     w_{j\ell} &\leq \rho^2_\ell\, y_j,  &&\forall \ell\in L, \; j \in J.
  \label{ctr:capacity2}
\end{align}
Constraints \eqref{ctr:flow2phaseB} enforce load conservation along the selected second-phase routes by linking the load carried by the large vehicle when leaving node $j \in J$ to that when leaving node $k \in J$. When arc $(j,k)$ is selected (i.e., $z_{jk\ell}=1$), these constraints reduce to $w_{k\ell} = w_{j\ell} + H_{k\ell}$, ensuring that the load at node $k$ equals the load at node $j \in J$ plus the amount of waste incorporated at node $k \in J$. 
When the arc is not selected, the terms involving $\rho^2_\ell$ relax the constraints.
Constraints \eqref{ctr:flow2phase0} enforce full unloading at the treatment plants, ensuring that the vehicle departs these nodes empty. Finally,  capacity constraints \eqref{ctr:capacity2}
limit the load leaving each consolidation facility to the vehicle’s capacity. They also ensure consistency with facility-locating decisions: if a consolidation facility is not located (i.e., $y_j = 0$), no load can originate from that node.

In the second phase, explicit connectivity constraints are not required. Because there is a single treatment plant for each waste type, the flow conservation constraints inherently ensure that all located consolidation facilities with positive inflow ($H_{j\ell} > 0$) are connected to the corresponding treatment plant to discharge the collected waste. This structural property guarantees network connectivity in the second phase without the need for additional constraints.

\subsection{Warm Start Solution}\label{sec:warm}

Finally, in this section, we present a procedure for constructing a feasible warm-start solution for the proposed MILP formulation for SWIFT. The motivation for developing such a procedure stems from the computational challenges associated with solving large-scale instances. In particular, realistic instances such as the Medell\'in case study considered in Section~\ref{sec:4} may generate very large branch-and-bound trees, making the discovery of good incumbent solutions an important factor in the efficiency of the solution process.

Providing an initial feasible solution is a well-established strategy for improving the performance of MILP solvers. Such warm-start solutions can accelerate the branch-and-bound procedure by supplying an incumbent bound from the beginning of the search, thereby facilitating node pruning and reducing computational effort. To this end, we propose a constructive warm-start procedure for SWIFT that exploits the structure of the problem. The procedure generates feasible solutions by decomposing the decision-making process into two main stages: the selection of consolidation facilities and treatment plants, and the construction of the collection and transfer routes in the two-echelon waste management network.

\paragraph{(i) Selection of consolidation facilities}
An initial set of consolidation facilities is selected using a \textit{greedy heuristic} based on average transportation costs. For each $i \in I$ and $j \in J$, we define the average transportation cost from collection point $i$ to consolidation facility $j$ as:
\[
\bar{c}_{ij} = \frac{1}{|L|} \sum_{\ell \in L} c^1_{ij\ell}.
\]
The first consolidation facility is selected as the one that minimizes the total average transportation cost from all collection points,
$j_1 = \arg\min_{j \in J} \sum_{i \in I} \bar{c}_{ij}.$

Additional consolidation facilities are then incorporated iteratively. Given a partial set of selected facilities $S \subseteq J$, the next facility is chosen as:
\[
j^* = \arg\min_{j \in J \setminus S}
\sum_{i \in I}
\min_{h \in S \cup \{j\}} \bar{c}_{ih}.
\]
The process continues while the installation budget remains feasible, namely,
\[
\sum_{j \in S} f_j^1
\leq
B - \sum_{\ell \in L} \max_{k \in K_\ell} \big\{f^2_{k\ell}\big\}.
\]

The selected set $S$ is then improved using a \textit{swap-based local search procedure}. At each iteration, a facility $j_{\text{out}} \in S$ is replaced by a facility $j_{\text{in}} \in J \setminus S$  whenever the resulting set $\hat{S} = \left(S \setminus \{j_{\text{out}}\}\right) \cup \{j_{\text{in}}\}$ leads to a lower total average transportation cost $\sum_{i \in I} \min_{h \in \hat{S}} \bar{c}_{ih}$. The procedure is repeated until no improving swap is found or a predefined maximum number of iterations is reached.

\paragraph{(ii) First-phase routing construction}
Given the set of open facilities $S$, for each waste type $\ell \in L$, each collection point $i \in I$ is  assigned to its nearest open facility according to:
\[
j(i,\ell) = \arg\min_{j \in S} c^1_{ij\ell}.
\]
This assignment induces, for each $j \in S$, a cluster of collection points that unload waste of type $\ell \in L$ at consolidation facility $j \in S$:
\[
C_j(\ell) = \{ i \in I : j(i,\ell) = j \}.
\]

A repair procedure is then applied to guarantee feasibility with respect to the second-phase vehicle capacity. In particular, if the total amount of waste type $\ell \in L$ assigned to a consolidation facility exceeds the corresponding second-phase vehicle capacity, collection points are iteratively reassigned to alternative facilities while preserving feasibility.

Then, routes within each cluster $C_j (\ell)$ are constructed using a \textit{nearest-neighbor heuristic}. Starting from consolidation facility $j \in J$, the first collection point visited by the small vehicle is selected as $i_1 = \arg\min_{i \in C_j(\ell)} c^1_{ji\ell}$.  Subsequently, the next collection point is determined as the closest unvisited point: 
\[
i_{k+1} =
\arg\min_{i \in C_j(\ell) \setminus \{i_1, \dots, i_k\}} c^1_{i_k i \ell}.
\]
Whenever the accumulated load exceeds the small-vehicle capacity, the route returns to the corresponding consolidation facility and a new subroute is initiated.

Additionally, the clusters are connected sequentially. Starting from a consolidation facility $j \in S$, the remaining open facilities are ordered according to a nearest-neighbor criterion. At each step, the next facility $j' \in S$ is selected among those not previously selected, minimizing the transportation cost $c^1_{jj'\ell}$ from the current facility. Then, for each pair of consecutive facilities $(j,j')$, a collection point $i^* \in C_{j'}(\ell)$ is selected to connect facility $j$ with cluster $C_{j'}(\ell)$, where $i^* = \arg\min_{i \in C_{j'}(\ell)} c^1_{ji\ell}$.

\paragraph{(iii) Second-phase routing construction}
Given the set of open consolidation facilities $S$, for each waste type $\ell \in L$, a treatment plant is selected from the set of candidate locations $K_\ell$ by minimizing the total transportation cost from facilities in $S$:
\[
k^*(\ell) = \arg\min_{k \in K_\ell} \sum_{j \in S} c^2_{jk\ell}.
\]
The consolidation facilities are then ordered using a \textit{nearest-neighbor heuristic} starting from the selected treatment plant. Specifically, at each step $t$, the next facility is selected as:
\[
j_{t+1} = \arg\min_{j \in S \setminus \{j_1,\dots,j_t\}} c^2_{j_tj\ell}.
\]
A route is subsequently constructed by visiting consolidation facilities sequentially. Whenever the accumulated load exceeds the large-vehicle capacity, the route returns to the treatment plant and a new trip is initiated.

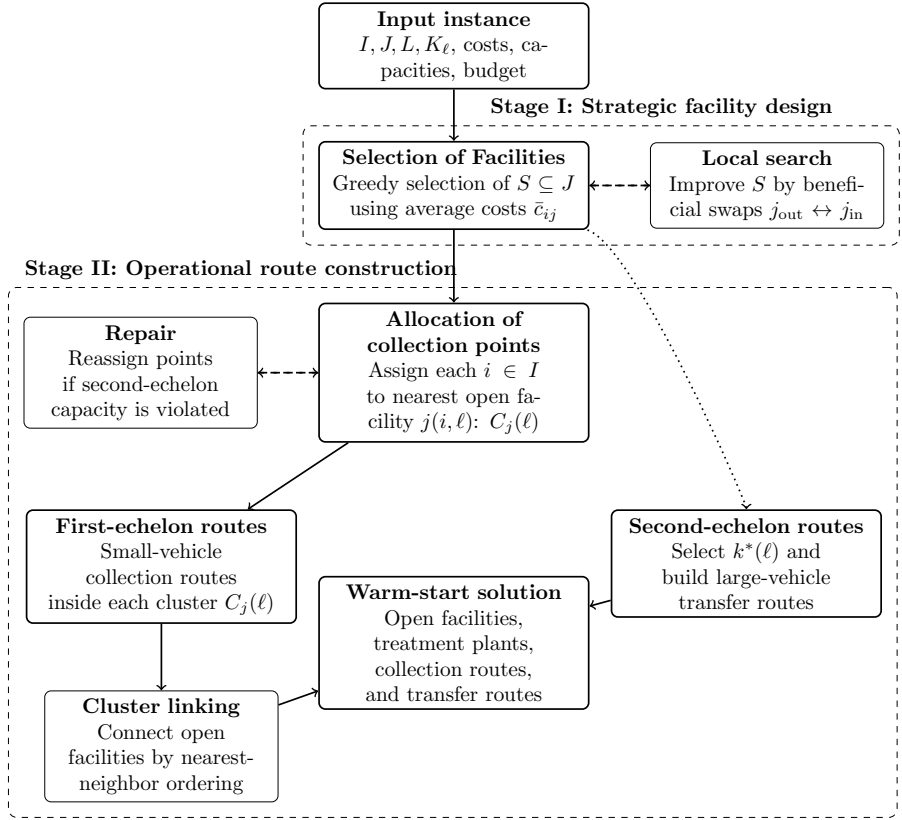
\begin{figure}[h]
\centering
\adjustbox{scale=0.80}{
\begin{tikzpicture}[
    font=\small,
    node distance=0.85cm and 1.0cm,
    main/.style={
        rectangle,
        rounded corners=3pt,
        draw,
        thick,
        align=center,
        text width=4.2cm,
        minimum height=1.05cm
    },
    smallbox/.style={
        rectangle,
        rounded corners=3pt,
        draw,
        align=center,
        text width=3.6cm,
        minimum height=0.9cm
    },
    phase/.style={
        rectangle,
        rounded corners=4pt,
        draw,
        dashed,
        inner sep=0.25cm
    },
    arrow/.style={->, thick},
    dashedarrow/.style={->, thick, dashed}
]

\node[main] (input) {
\textbf{Input instance}\\
$I,J,L, K_\ell$, costs, capacities, budget
};

\node[main, below=of input] (facilities) {
\textbf{Selection of Facilities }\\
Greedy selection of $S\subseteq J$ using average costs
$\bar c_{ij}$
};

\node[smallbox, right=of facilities] (swap) {
\textbf{Local search}\\
Improve $S$ by beneficial swaps
$j_{\rm out}\leftrightarrow j_{\rm in}$
};

\node[main, below=1.2cm of facilities] (assign) {
\textbf{Allocation of\\ collection points }\\
Assign each $i\in I$ to nearest open facility
$j(i,\ell)$: $C_j(\ell)$
};

\node[smallbox, left=of assign] (repair) {
\textbf{Repair}\\
Reassign points if second-echelon \\ capacity is violated
};

\node[main, below left=1.1cm and 0.35cm of assign] (phase1) {
\textbf{First-echelon routes}\\
Small-vehicle \\collection routes\\ inside each cluster $C_j(\ell)$
};

\node[main, below right=1.1cm and 0.35cm of assign] (phase2) {
\textbf{Second-echelon routes}\\
Select $k^*(\ell)$ and build large-vehicle transfer routes
};

\node[smallbox, below=1.1cm of phase1] (connect) {
\textbf{Cluster linking}\\
Connect open\\ facilities by nearest-neighbor ordering
};

\node[main, below=2.2cm of assign] (output) {
\textbf{Warm-start solution}\\
Open  facilities, \\treatment plants, collection routes, and transfer routes
};

\draw[arrow] (input) -- (facilities);
\draw[arrow] (facilities) -- (assign);
\draw[arrow] (assign) -- (phase1);
\draw[arrow,dotted]
    (facilities) .. controls +(3,-1) and +(0,1) .. (phase2);
\draw[arrow] (phase1) -- (connect);
\draw[arrow] (connect) -- (output);
\draw[arrow] (phase2) -- (output);

\draw[dashedarrow] (facilities) -- (swap);
\draw[dashedarrow] (swap) -- (facilities);

\draw[dashedarrow] (assign) -- (repair);
\draw[dashedarrow] (repair) -- (assign);

\node[
    phase,
    fit=(facilities)(swap),
    label={[xshift=1cm]above:\textbf{Stage I: Strategic facility design}}
] {};
\node[
phase, 
fit=(assign)(repair)(phase1)(phase2)(connect),
label={[xshift=-3.5cm]above:\textbf{Stage II: Operational route construction}}
] {};

\end{tikzpicture}}
\caption{Structure of the constructive warm-start algorithm for SWIFT.}
\label{fig:warm_start_flowchart}
\end{figure}

Overall, the proposed heuristic guarantees feasibility with respect to the capacity constraints and provides a feasible warm-start solution for SWIFT. A schematic overview of the procedure is presented in the flowchart of Figure \ref{fig:warm_start_flowchart}.

\section{Case Study: Urban Waste Management in Medell\'in}\label{sec:4}

This section evaluates SWIFT through a real-world case study based on the urban solid waste management system of Medellín, Colombia. The setting involves multiple recyclable waste streams that must be collected across the city, unloaded at intermediate consolidation facilities, and transported to specialized treatment plants.

The case study builds on the Medellín instance described in \citet{rodriguez2026bi}, which focuses on strategic planning decisions in urban waste management. Medellín provides a suitable test bed due to its dense urban structure, heterogeneous spatial distribution of waste generation, and the need to coordinate small-scale collection with larger transfer operations. In contrast to purely strategic approaches, the proposed SWIFT integrates facility location, routing decisions, vehicle capacity limitations, and waste-type-specific transportation requirements within a unified framework.

\subsection{Experimental Settings}
\label{subsec:medellin_dataset}

The instance consists of $n=471$ collection points, $m=16$ candidate consolidation facilities,   $L=4$  types of waste, and $k_\ell$ candidate treatment plants (or depots) for each waste type, with $k_\ell = 21, 12, 3, 5$. Collection points represent waste generation locations distributed throughout the city, while candidate consolidation facilities correspond to intermediate sites where small vehicles can unload collected waste before transferring it to treatment plants.

\begin{figure}[h]
\begin{center}
\includegraphics[width=0.7\textwidth]{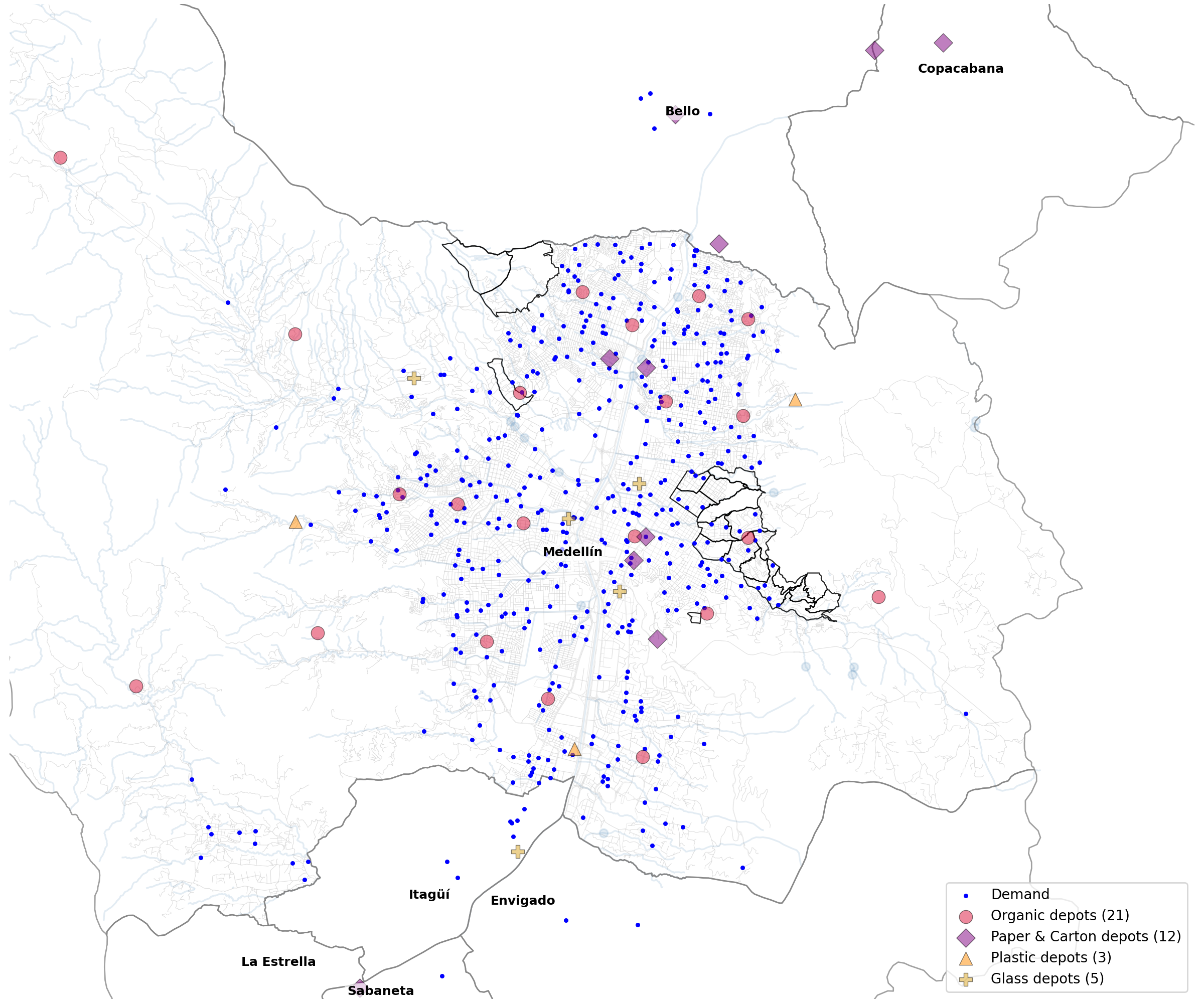}
\end{center}
\caption{Geographical distribution of the facilities in the Medell\'in area used in our experiments.\label{fig:map}}
\end{figure}

The geographical distribution of the different types of facilities is shown in Figure~\ref{fig:map}.

The four types of waste correspond to the main recyclable or valorizable streams: Organic, Paper and Carton, Plastic, and Glass. For each type, waste is collected by a dedicated small vehicle and subsequently transported by a larger vehicle from the selected consolidation facilities to a treatment plant. This two-level structure captures the operational distinction between collection in dense urban areas and long-haul transportation to specialized processing sites.

The main goal of this study is to assess how the integrated optimization framework supports the design of efficient recycling networks by jointly determining facility locations and routing decisions. In particular, we analyze the trade-offs between infrastructure investment and transportation efficiency under realistic operational constraints.

To this end, the model is solved for four infrastructure configurations, $p\in\{8,10,12,14\}$, where $p$ denotes the maximum number of consolidation facilities that can be located, and then, representing the investment to be done in $p$ similar types of consolidation points (see \eqref{mC:ctr2b}). These values range from moderately centralized to highly distributed infrastructure layouts, enabling the evaluation of how consolidation intensity affects routing efficiency, operational costs, and facility utilization.

All computational experiments were performed on Huawei FusionServer Pro XH 321 \texttt{albaic\'in} at Universidad de Granada (\url{supercomputacion.ugr.es/arquitecturas/albaicin}) with an Intel Xeon Gold 6258R CPU @ 2.70GHz with 28 cores. Optimization tasks were solved using Gurobi Optimizer version~12.0.1 within a time limit of 12 hours.

\subsection{Obtained Results}
\label{subsec:baseline_results}

This subsection analyzes the solutions obtained for the baseline Medellín dataset under the proposed SWIFT framework. The analysis focuses on three complementary aspects of the obtained solutions:
(i) the spatial structure of the integrated location-routing design,
(ii) the operational efficiency generated by additional consolidation infrastructure,
(iii) the robustness of the resulting network under alternative urban operating scenarios, and (iv) sustainability insights of the obtained solutions.

The results provide insights not only into the routing behavior induced by the proposed optimization framework, but also into the role of intermediate consolidation infrastructure in improving transportation efficiency, reducing urban collection effort, and enhancing the operational sustainability of multi-product recycling systems.

\subsubsection{Spatial Structure of the Obtained Solutions}

\begin{figure}[t]
\centering
\includegraphics[width=0.95\textwidth]{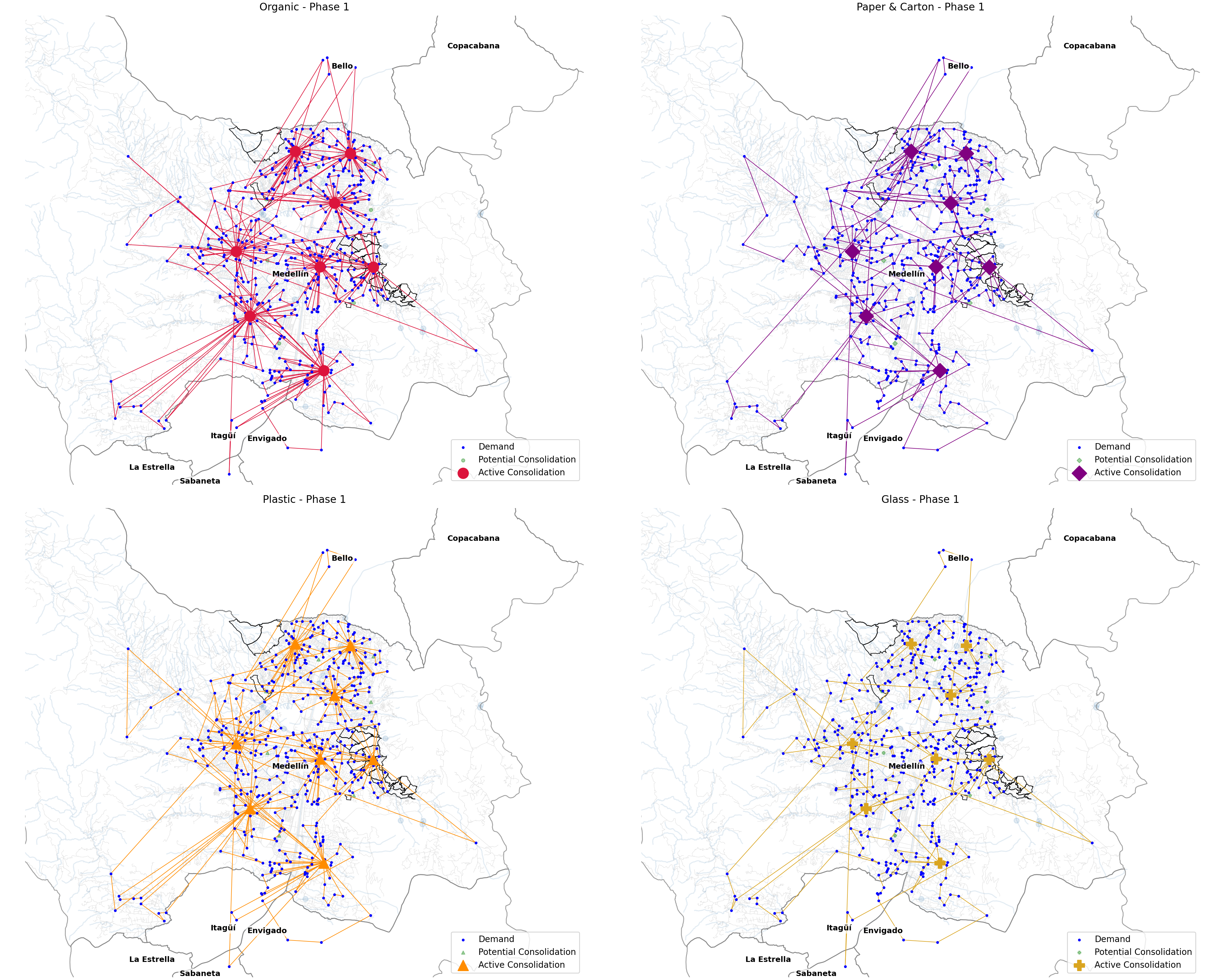}
\caption{Collection-phase routing structure for the baseline Medellín dataset with $p=8$ activated consolidation facilities.}
\label{fig:phase1_8}
\end{figure}

\begin{figure}[t]
\centering
\includegraphics[width=0.95\textwidth]{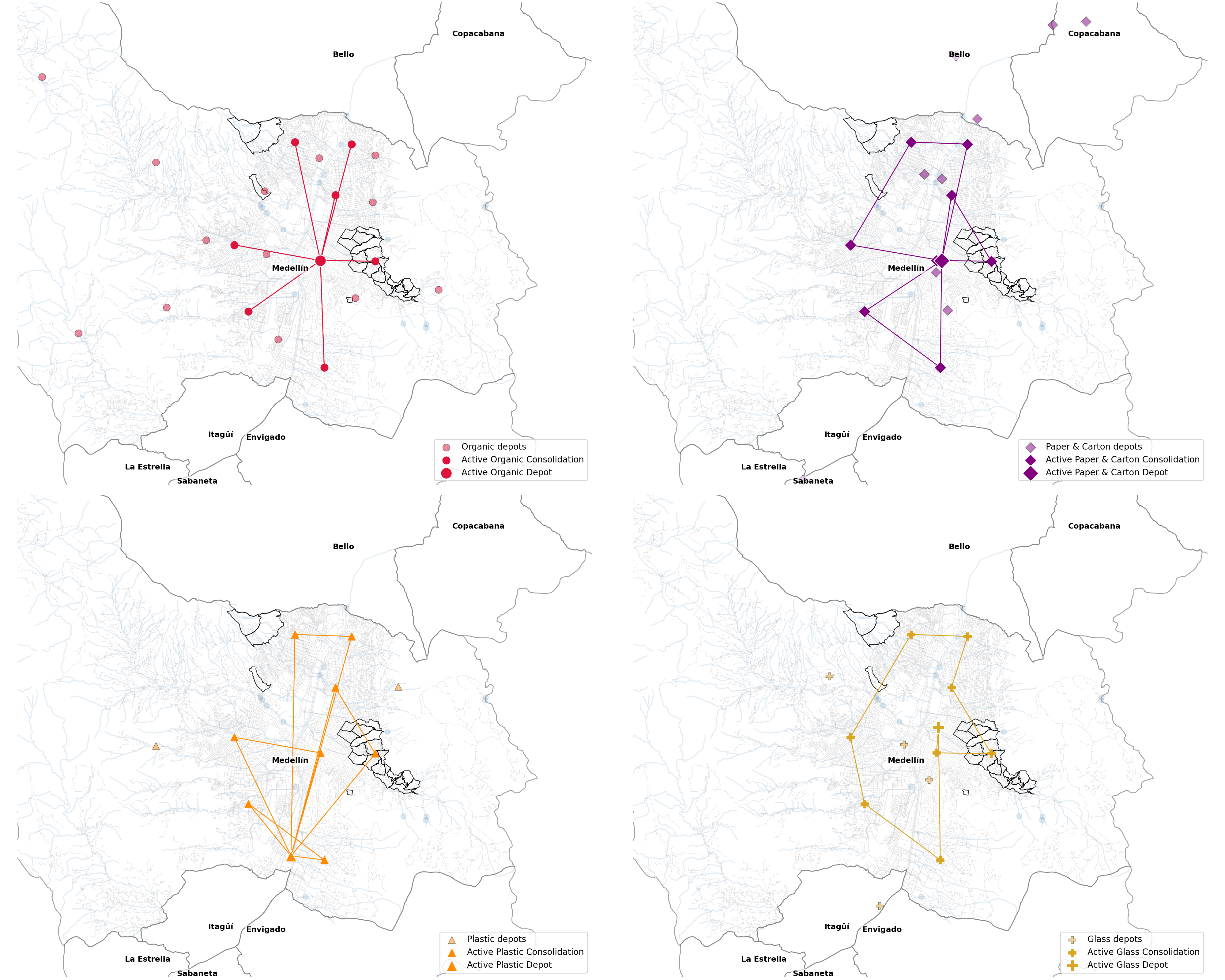}
\caption{Transfer-phase routing structure for the baseline Medellín dataset with $p=8$ activated consolidation facilities.}
\label{fig:phase2_8}
\end{figure}

We first analyze the baseline configuration with $p=8$ activated consolidation facilities, which serves as the reference benchmark for the remaining experiments. The associated routing structures for the collection and transfer phases are illustrated in Figures~\ref{fig:phase1_8} and~\ref{fig:phase2_8}, respectively.

Each waste stream is represented independently using a specific combination of color and symbol: red circles for Organic waste, purple diamonds for Paper \& Carton, orange triangles for Plastic, and yellow crosses for Glass. Activated consolidation facilities are represented with solid medium-sized symbols, candidate but non-activated facilities appear as transparent symbols, and activated treatment depots are highlighted using larger symbols. The displayed arcs represent the optimized vehicle routes obtained for each waste stream and transportation phase.

The first important observation concerns the spatial distribution of the activated consolidation facilities. The selected locations are geographically dispersed across the urban area, generating a balanced infrastructure capable of efficiently supporting both transportation phases. These facilities simultaneously play two interconnected roles: they operate as unloading hubs for small collection vehicles during the first transportation phase and as aggregation nodes for the transfer operations toward treatment depots during the second phase.

Consequently, the selected locations emerge from a nontrivial trade-off between local collection efficiency and long-haul transportation efficiency. On the one hand, consolidation facilities should be sufficiently close to the demand points to reduce the travel distance of small collection vehicles operating under tight capacity restrictions. On the other hand, they should also remain well positioned with respect to the treatment depots in order to limit the cost of the transfer phase. The resulting infrastructure configuration therefore reflects the integrated nature of the proposed optimization framework.

The location of the treatment depots also exhibits meaningful operational patterns. Since Organic waste generates the largest transportation volumes and the most intensive routing activity, its associated depot is positioned close to the central demand region, thereby reducing the cost of heavily utilized transfer routes. Paper \& Carton and Glass depots also exhibit relatively central locations. In contrast, the limited availability of feasible candidate locations for Plastic treatment facilities forces the model to place its depot in the southern region of Medellín, producing a noticeably different transfer structure.

Figure~\ref{fig:phase1_8} additionally reveals a clear spatial clustering behavior around the activated consolidation facilities. These facilities naturally induce compact collection regions, generating relatively short and geographically coherent routes. Such behavior is particularly advantageous under the limited capacities of the small collection vehicles, since it reduces travel effort, unloading frequency, and urban congestion exposure.

The routing structures also differ substantially across waste streams due to the interaction between demand intensity, vehicle capacities, and facility locations. Organic waste generates highly localized routes in which vehicles typically visit only a small number of demand points before unloading. This behavior is consistent with the relatively large accumulated loads associated with this stream.

Paper \& Carton and Plastic exhibit intermediate routing structures. Their lower demand intensity allows vehicles to visit a larger number of customers before unloading, thereby generating routes with wider geographical coverage while preserving coherent service regions around the activated consolidation facilities.

Glass exhibits the most spatially extended collection routes. Since transported volumes are substantially smaller, vehicles are capable of visiting many demand points before reaching their capacity limits. This generates elongated and interconnected collection tours across the urban area.

Overall, these results illustrate the ability of the proposed SWIFT formulation to simultaneously adapt facility-location and routing decisions to the operational characteristics of each waste stream. Rather than imposing predefined routing structures, the model endogenously determines the most appropriate collection patterns according to vehicle utilization, spatial demand distribution, and infrastructure availability.

The second transportation phase, illustrated in Figure~\ref{fig:phase2_8}, presents a substantially different operational structure. In this phase, large-capacity vehicles transport aggregated waste flows from the consolidation facilities toward the corresponding treatment depots. Since the material has already been consolidated, the resulting routes involve fewer visited nodes and exhibit significantly simpler topological patterns.

The transfer structures differ considerably across waste streams. Organic waste generates an almost star-shaped configuration composed primarily of direct trips between consolidation facilities and depots. This indicates that vehicles typically become fully loaded after servicing a single consolidation facility, making multi-stop transfer routes operationally inefficient.

Paper \& Carton exhibits an intermediate behavior characterized by multiple cyclic routes connecting subsets of consolidation facilities before returning to the treatment depot. This structure balances transportation efficiency and vehicle-capacity utilization while preserving spatial coherence.

Plastic displays a slightly more fragmented configuration, requiring additional transfer trips due to the combined effect of spatial dispersion and accumulated transported volume.

In contrast, Glass produces the most consolidated transfer structures. In several cases, a single transfer vehicle is capable of visiting multiple consolidation facilities within the same route before unloading at the depot. This behavior reflects the comparatively low transported volumes associated with this waste stream and generates highly efficient cyclic transportation patterns.

From a sustainability perspective, these results reveal the importance of jointly optimizing collection and transfer operations. The proposed framework coordinates routing structures with the operational characteristics of each waste stream, thereby reducing unnecessary travel distances and improving vehicle utilization. Such coordinated decisions are particularly relevant in dense urban environments, where transportation efficiency directly affects fuel consumption, emissions, congestion, and operational costs.

\subsubsection{Infrastructure Efficiency and Marginal Benefits}

An important managerial insight concerns the role of the activated consolidation infrastructure. Table~\ref{tab:improvement_original} reports the relative reduction in first-phase collection costs obtained by increasing the number of activated consolidation facilities with respect to the baseline configuration with $p=8$. In this baseline scenario, the transportation cost associated with the first phase (i.e., a single-day collection operation covering all waste demands) was \$10,790,253 Colombian pesos (approximately \$2,250 USD), whereas the transportation cost of the second phase (i.e., a single-day operation in which large vehicles collect waste from the consolidation facilities and transport it to the corresponding treatment plants) was \$872,027 Colombian pesos (approximately \$182 USD).
\begin{table}[h]
\centering
\renewcommand{\arraystretch}{1.15}
\small
\begin{tabular}{cccc}
\toprule
$p=8$ & $p=10$ & $p=12$ & $p=14$ \\
\midrule
0.00\% & 17.73\% & 19.89\% & 20.27\% \\
\bottomrule
\end{tabular}
\caption{Relative improvement in first-phase collection costs with respect to the baseline configuration with $p=8$ activated consolidation facilities. \label{tab:improvement_original}}
\end{table}
The results reveal a substantial reduction in collection costs as additional consolidation facilities are installed. Increasing the number of facilities from $p=8$ to $p=10$ reduces first-phase routing costs by approximately $17.7\%$. Further increases to $p=12$ and $p=14$ produce improvements of approximately $19.9\%$ and $20.3\%$, respectively.

These reductions are primarily explained by the improved spatial accessibility generated by the additional infrastructure. A denser consolidation network allows collection vehicles to unload closer to the collection points, thereby reducing collection-route lengths and improving operational efficiency under strict vehicle-capacity restrictions.

At the same time, the results clearly exhibit diminishing marginal returns. While the transition from $p=8$ to $p=10$ generates substantial efficiency gains, the additional improvements obtained beyond $p=12$ become comparatively moderate. This behavior indicates that, after a certain infrastructure density threshold, the collection regions already become sufficiently compact and additional facilities contribute only limited routing improvements.

This observation is particularly relevant from a sustainability and investment-planning perspective. Although additional consolidation facilities improve operational efficiency, they also require additional infrastructure investment, land usage, and operational maintenance. The results therefore suggest the existence of a balanced infrastructure region in which the marginal operational benefits approximately stabilize.

\begin{figure}[t]
\centering
\includegraphics[width=0.95\textwidth]{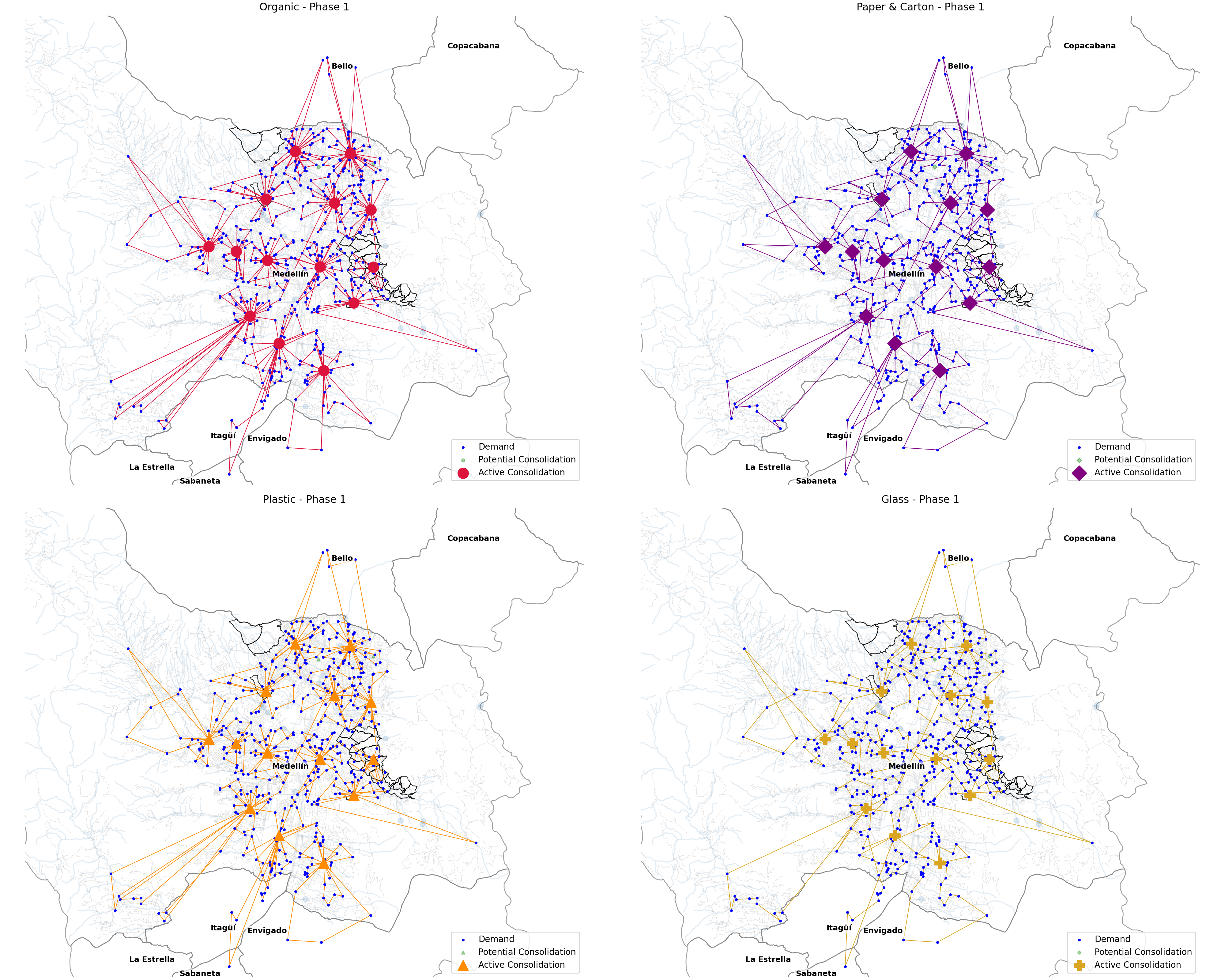}
\caption{Collection-phase routing structure for the baseline dataset with $p=14$ activated consolidation facilities.}
\label{fig:phase1_14}
\end{figure}

\begin{figure}[t]
\centering
\includegraphics[width=0.95\textwidth]{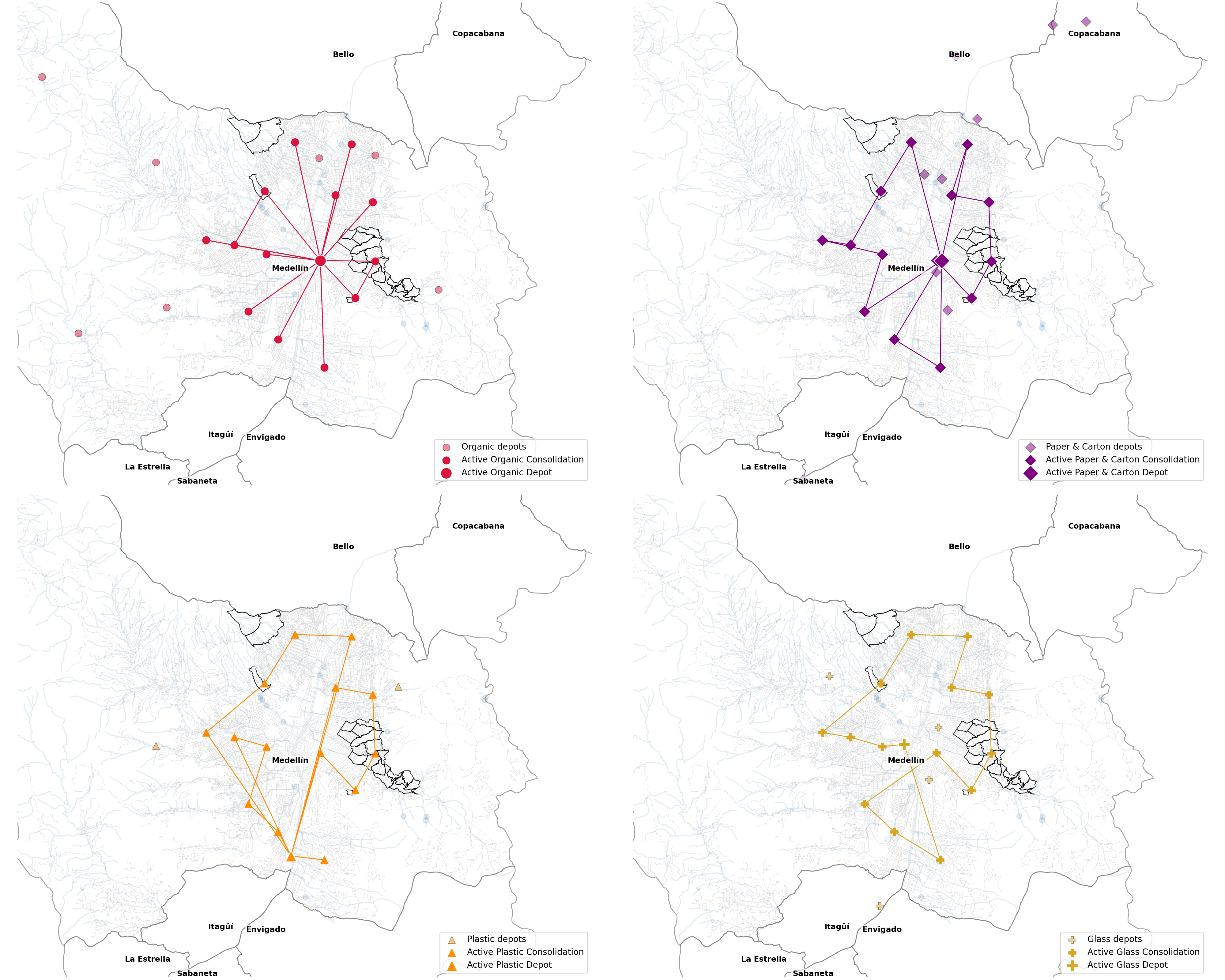}
\caption{Transfer-phase routing structure for the baseline dataset with $p=14$ activated consolidation facilities.}
\label{fig:phase2_14}
\end{figure}

Figures~\ref{fig:phase1_14} and~\ref{fig:phase2_14} illustrate the routing structures obtained for the densest tested infrastructure configuration ($p=14$). 

Figures~\ref{fig:phase1_8} and~\ref{fig:phase1_14} allow a direct comparison of the collection-phase routing structures obtained under two different infrastructure levels. When the number of activated consolidation facilities increases from $p=8$ to $p=14$, the collection network undergoes a clear spatial reorganization, resulting in smaller and more compact service regions.

For $p=8$, each consolidation facility is responsible for a relatively large geographical area. Consequently, collection routes are longer and exhibit a more radial structure, particularly for Organic and Paper \& Carton waste. Several routes connect distant demand points with the same consolidation facility, generating elongated collection patterns and substantial overlap between service regions. This behavior is especially visible in peripheral areas of the city, where vehicles must travel long distances before unloading.

In contrast, the configuration with $p=14$ produces a significantly denser infrastructure network. The additional consolidation facilities allow the model to partition the urban area into smaller and more balanced service regions. As a result, collection routes become substantially shorter and geographically concentrated around nearby unloading points. This effect is particularly evident for Organic waste, where the routes evolve from large radial structures into highly localized clusters centered around the activated facilities.

The impact of the additional infrastructure is also noticeable for the remaining waste streams. For Paper \& Carton and Plastic, the denser consolidation network reduces the need for long cross-city connections and generates more coherent local routing structures. In the case of Glass, although the routes remain spatially more extensive due to the lower accumulated loads, the additional facilities still reduce the average collection distance and improve route compactness.

From an operational perspective, these results indicate that increasing the number of consolidation facilities substantially improves the efficiency of the first transportation phase by reducing collection distances, unloading trips, and vehicle exposure to urban congestion. At the same time, the comparison highlights the existence of diminishing structural changes beyond a certain infrastructure density. While the transition from $p=8$ to $p=14$ clearly improves route compactness and territorial balance, several service regions already appear reasonably stabilized for intermediate infrastructure levels, suggesting that additional facilities eventually provide progressively smaller operational gains.

Figures~\ref{fig:phase2_8} and~\ref{fig:phase2_14} provide a comparative visualization of the transfer-phase routing structures obtained for two different infrastructure levels. In contrast to the collection phase, where the main impact of increasing the number of consolidation facilities was the creation of more compact local service regions, the second transportation phase exhibits a more heterogeneous structural evolution across waste streams.

For the baseline configuration with $p=8$, the transfer network is relatively sparse and highly centralized. In particular, Organic waste displays an almost pure star-shaped configuration, where each consolidation facility is connected directly to the treatment depot through an independent transfer route. This behavior indicates that the accumulated waste volumes at the consolidation facilities are sufficiently large to saturate the transfer vehicles individually, making route consolidation operationally unattractive.

When the number of consolidation facilities increases to $p=14$, the transfer structure becomes more interconnected and spatially distributed. Since the collected waste is now divided among a larger number of consolidation facilities, the amount of material accumulated at each facility decreases. Consequently, transfer vehicles can combine multiple facilities within the same route before unloading at the depot. This effect is especially visible for Paper \& Carton and Glass, where the routes evolve from partially radial structures into more cyclic and interconnected configurations.

Plastic exhibits an intermediate behavior. Although the denser infrastructure reduces the average amount of waste aggregated at each consolidation point, the geographical dispersion of the facilities still induces relatively elongated transfer movements. Nevertheless, the routing structure for $p=14$ remains more spatially balanced and operationally flexible than the corresponding configuration with $p=8$.

The most notable structural transformation occurs for Glass. Under $p=8$, several transfer routes already combine multiple consolidation facilities, but the larger infrastructure density associated with $p=14$ allows the formation of even more consolidated cyclic routes. This indicates that the lower transported volumes associated with Glass strongly favor route integration during the transfer phase.

Overall, the comparison between Figures~\ref{fig:phase2_8} and~\ref{fig:phase2_14} highlights an important operational trade-off induced by additional consolidation infrastructure. Increasing the number of facilities substantially improves the efficiency of the collection phase by reducing local transportation distances, but it also generates a larger number of transfer nodes that must be coordinated during the second phase. Consequently, the transfer network becomes more interconnected and operationally complex as the infrastructure density increases.

Taken together, these results highlight the benefits of the integrated optimization framework. By jointly optimizing facility-location and routing decisions across both transportation phases, the proposed model generates spatially coherent collection regions, efficient transfer structures, and waste-specific transportation patterns that would be difficult to obtain through sequential planning approaches.

\subsubsection{Sensitivity Analysis}
\label{subsec:scenario_design}

To evaluate the robustness of the obtained network designs and derive additional managerial insights, we constructed a set of representative operating scenarios obtained by perturbing selected parameters of the baseline Medellín instance, denoted as Original.

The network topology, candidate locations, and compatibility relations remain fixed across all experiments, so that each scenario represents an alternative operational environment characterized by different demand, capacity, or transportation conditions.

\begin{table}[t]
\centering
\renewcommand{\arraystretch}{1.05}
\adjustbox{scale=0.60}{
\begin{tabular}{lcccc}
\toprule
\textbf{Scenario} 
& \textbf{Demand $d_{i\ell}$} 
& \textbf{Spatial distribution} 
& \textbf{Vehicle capacities} 
& \textbf{Transportation costs} \\
\midrule
Original 
& baseline 
& baseline 
& baseline 
& baseline \\

\textit{High Volume} 
& $1.25\,d_{i\ell}$ 
& unchanged 
& baseline 
& baseline \\

\textit{Skewed Mix} 
& $1.5\,d_{i\ell^\star}$ 
& unchanged 
& baseline 
& baseline \\

\textit{Urban Concentration} 
& total preserved 
& concentrated demand 
& baseline 
& baseline \\

\textit{Capacity Stress} 
& baseline 
& unchanged 
& $0.85\,(\rho^1_\ell,\rho^2_\ell)$ 
& baseline \\

\textit{Cost Inflation} 
& baseline 
& unchanged 
& baseline 
& $1.2\,(c^1_{ij\ell},c^2_{jk\ell})$ \\
\bottomrule
\end{tabular}}
\caption{Scenario definition and parameter perturbations relative to the baseline Medellín instance.}
\label{tab:scenarios_medellin}
\end{table}

The \textit{High Volume} scenario represents increased waste generation associated with population growth or consumption changes. The \textit{Skewed Mix} scenario captures changes in waste composition in which one material stream becomes dominant. The \textit{Urban Concentration} scenario models spatial concentration of waste generation in dense urban areas. The \textit{Capacity Stress} scenario reduces vehicle capacities, while the \textit{Cost Inflation} scenario represents generalized transportation-cost increases.



\begin{figure}[t]
\centering
\includegraphics[width=0.95\textwidth]{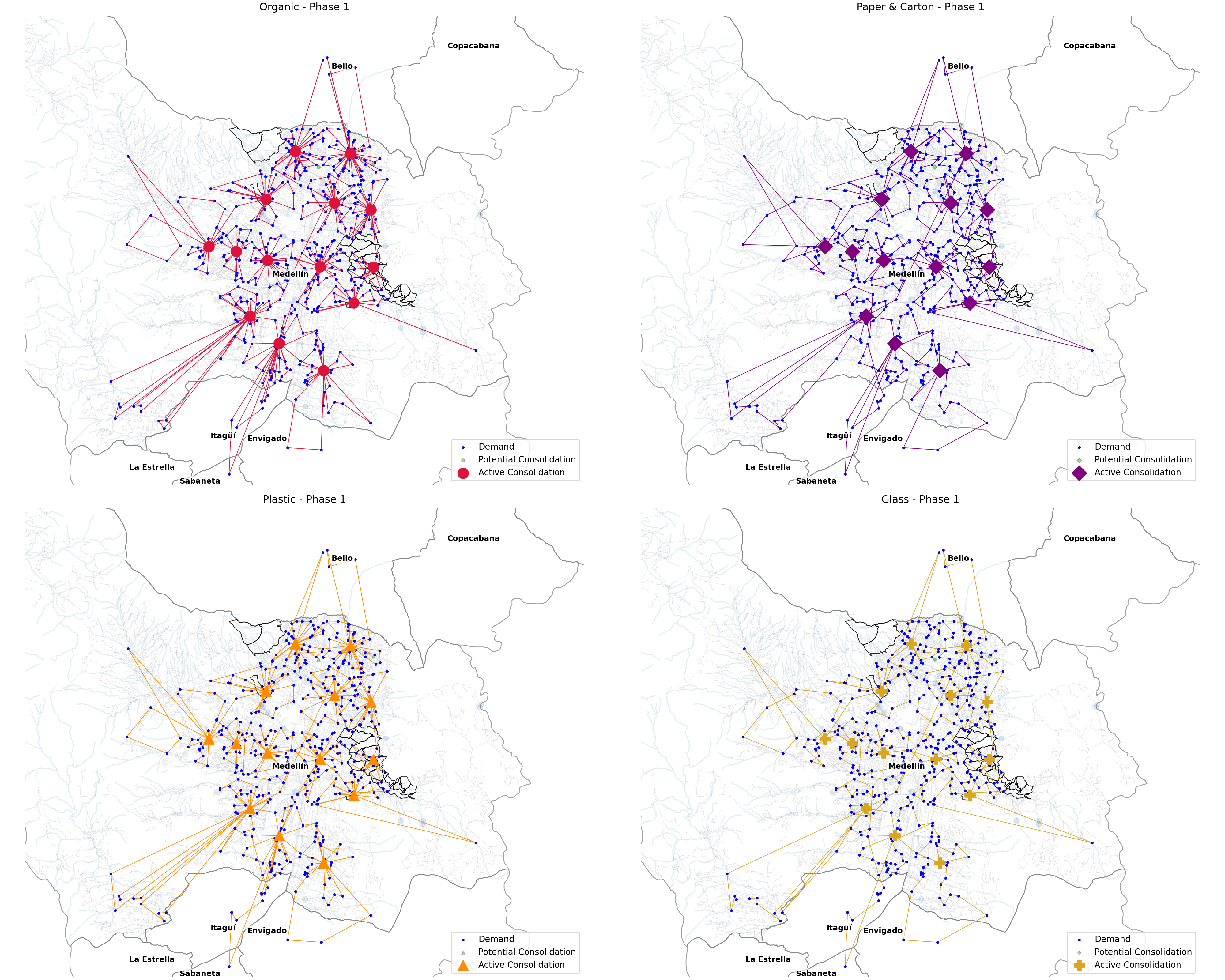}
\caption{Collection-phase routing structure for the baseline dataset with $p=14$ activated consolidation facilities  for the UrbanConcentration scenario..}
\label{fig:phase1_14_urban}
\end{figure}

\begin{figure}[t]
\centering
\includegraphics[width=0.95\textwidth]{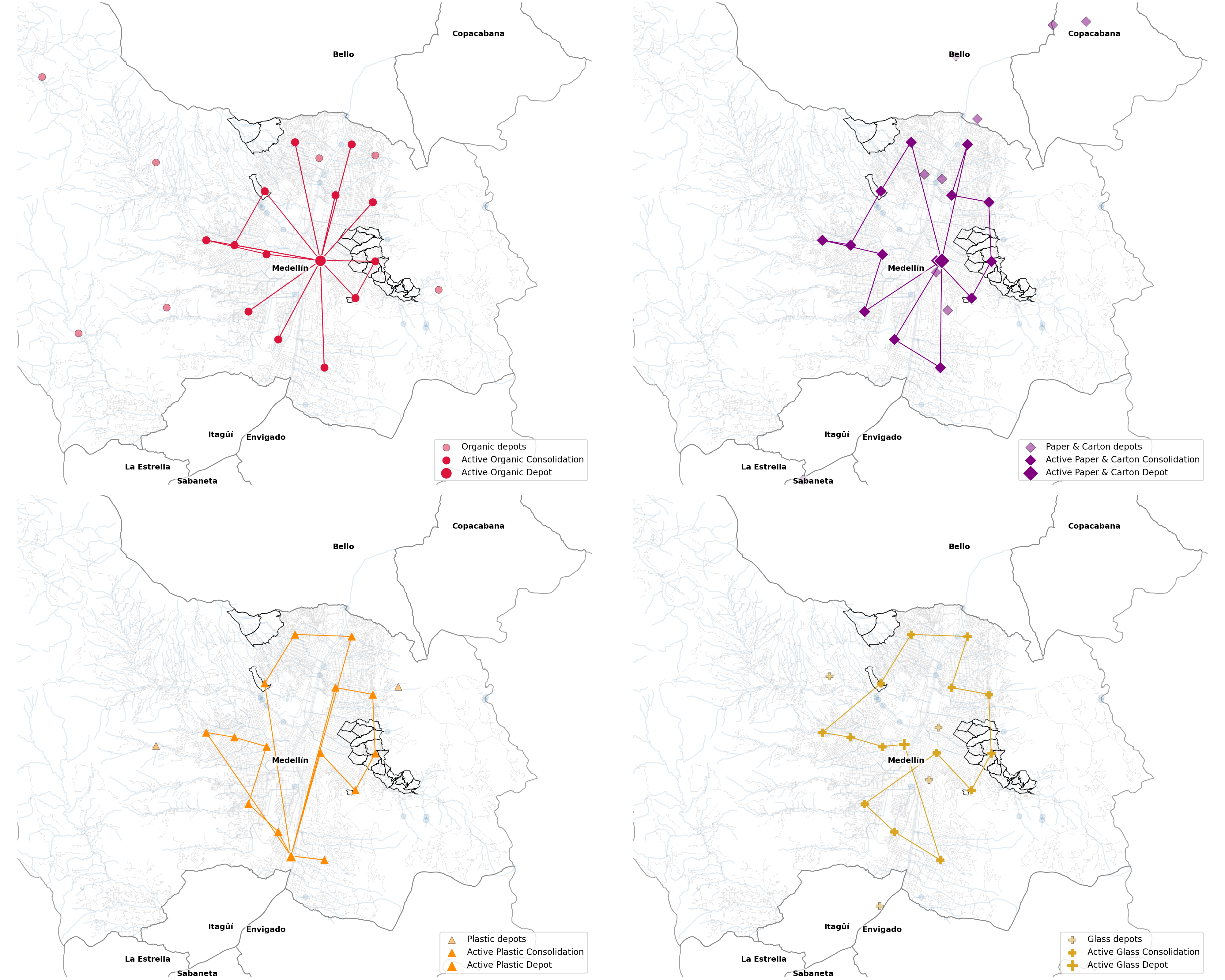}
\caption{Transfer-phase routing structure for the baseline dataset with $p=14$ activated consolidation facilities for the UrbanConcentration scenario.}
\label{fig:phase2_14_urban}
\end{figure}

In Figures \ref{fig:phase1_14_urban} and \ref{fig:phase2_14_urban} we show the graphical solutions for the scenario UrbanConcentration. Under this scenario, the demand points become more spatially clustered around the central urban area of Medell\'in. As a consequence, the routing structures generated by SWIFT become noticeably denser and geographically more compact than in the Original scenario. In particular, the routes associated with Organic and Paper \& Carton waste exhibit shorter radial extensions and stronger overlap around the central districts, reflecting the higher spatial proximity between demand points and consolidation facilities.

The reduction in spatial dispersion also affects the operational organization of the collection phase. In the Original scenario, several routes extend toward peripheral regions of the city, requiring long collection trips before unloading. In contrast, the \textit{Urban Concentration} scenario generates more localized service regions, reducing the need for long-distance urban movements and increasing route compactness across all waste streams.

This effect is especially visible for Plastic and Glass, where the collection tours become substantially more interconnected within the urban core. Since these waste streams are associated with lower effective loads, the increased spatial concentration allows vehicles to visit multiple nearby demand points with reduced travel effort before reaching capacity limits.

From a managerial perspective, these results indicate that urban-density patterns strongly influence the operational efficiency of multi-product recycling systems. Spatially concentrated waste generation naturally favors compact collection structures and reduces the transportation burden associated with the first phase of the system.



Table~\ref{tab:phase1_share} reports the relative contribution of the collection and transfer phases to the total operational cost across the different scenarios and infrastructure configurations.

\begin{table}[t]
\centering

\renewcommand{\arraystretch}{1.10}
\adjustbox{scale=0.75}{
\begin{tabular}{lrrr|lrrr}
\toprule
\textbf{Scenario} & \textbf{$p$} & \textbf{Phase 1} & \textbf{Phase 2} 
& \textbf{Scenario} & \textbf{$p$} & \textbf{Phase 1} & \textbf{Phase 2} \\
\midrule
Original 
& 8  & 92.52\% & 7.48\%
& \textit{Cost Inflation} 
& 8  & 92.33\% & 7.67\% \\
& 10 & 90.63\% & 9.37\%
&  
& 10 & 90.67\% & 9.33\% \\
& 12 & 89.46\% & 10.54\%
&  
& 12 & 89.46\% & 10.54\% \\
& 14 & 89.05\% & 10.95\%
&  
& 14 & 89.05\% & 10.95\% \\
\midrule
\textit{High Volume} 
& 8  & -- & --
& \textit{Skewed Mix} 
& 8  & -- & -- \\
& 10 & 92.16\% & 7.84\%
&  
& 10 & -- & -- \\
& 12 & 89.24\% & 10.76\%
&  
& 12 & 91.00\% & 9.00\% \\
& 14 & 88.57\% & 11.43\%
&  
& 14 & 89.60\% & 10.40\% \\
\midrule
\textit{Capacity Stress} 
& 8  & -- & --
& \textit{Urban Concentration} 
& 8  & 91.42\% & 8.58\% \\
& 10 & 91.20\% & 8.80\%
&  
& 10 & 89.56\% & 10.44\% \\
& 12 & 89.47\% & 10.53\%
& 
& 12 & 89.37\% & 10.63\% \\
& 14 & 89.14\% & 10.86\%
&  
& 14 & 89.24\% & 10.76\% \\
\bottomrule
\end{tabular}}
\caption{Relative contribution of collection and transfer phases to the total operational cost across scenarios and infrastructure configurations. \label{tab:phase1_share}}
\end{table}

Entries marked with ``--'' correspond to infeasible instances, indicating that the available consolidation infrastructure is insufficient to support the corresponding demand and vehicle-capacity conditions. These infeasibilities appear primarily in the \textit{High Volume}, \textit{Capacity Stress}, and \textit{Skewed Mix} scenarios for small values of $p$, revealing the existence of a minimum infrastructure threshold required to maintain operational feasibility under stressed urban conditions.

For all feasible instances, the collection phase systematically dominates the operational cost, accounting for approximately $89\%$-$93\%$ of the total objective value. This behavior reflects the intrinsic complexity of dense urban collection operations, where small vehicles must visit all demand points under heterogeneous waste-generation patterns and capacity limitations.

An important structural trend emerges as the number of consolidation facilities increases. The relative contribution of the collection phase decreases consistently with larger values of $p$, whereas the contribution of the transfer phase increases. This behavior reflects the trade-off generated by additional infrastructure. A denser consolidation network reduces local collection distances and improves unloading accessibility, but simultaneously increases the number of transfer operations connecting consolidation facilities with treatment depots.

The observed pattern remains remarkably stable across all tested scenarios, indicating that the balance between collection and transfer operations is primarily driven by the intrinsic hierarchical structure of the recycling system.

\subsubsection{Sustainability}

From a sustainability perspective, the results indicate that intermediate consolidation infrastructure can play a fundamental role in improving the environmental and operational performance of urban recycling systems. Increasing the number of activated consolidation facilities generates more compact collection regions and substantially reduces the distance traveled by small collection vehicles during the first transportation phase. This effect is particularly relevant in dense urban environments, where fuel consumption, emissions, noise, and congestion externalities are typically more severe. Consequently, the proposed SWIFT framework contributes not only to reducing operational costs, but also to mitigating the environmental footprint associated with municipal waste collection activities.

The experiments additionally show that the benefits of a denser infrastructure network extend beyond the collection phase. By distributing the accumulated waste among a larger number of consolidation facilities, the transfer phase gains additional operational flexibility, allowing transfer vehicles to combine multiple facilities within the same route and improving vehicle utilization in several waste streams. In particular, Paper \& Carton and Glass exhibit more consolidated and interconnected transfer structures as the infrastructure density increases, reducing the prevalence of partially loaded trips and improving transportation efficiency.

Nevertheless, the results also reveal the existence of important trade-offs associated with infrastructure expansion. Although additional consolidation facilities improve route compactness and reduce collection effort, they simultaneously increase the complexity of the transfer network by creating additional inter-facility transportation requirements. Moreover, infrastructure expansion entails additional investment costs, land usage, operational maintenance, and potential urban-space impacts. The observed diminishing marginal gains obtained for larger values of $p$ suggest that, beyond a certain infrastructure density, additional facilities provide progressively smaller operational improvements.

These observations highlight the importance of jointly optimizing infrastructure and routing decisions within an integrated framework. Rather than simply minimizing transportation costs, the proposed SWIFT model supports the identification of balanced logistics configurations capable of simultaneously improving collection efficiency, vehicle utilization, operational feasibility, and environmental sustainability. Such integrated planning tools are particularly valuable for supporting long-term sustainable urban waste-management policies in rapidly growing metropolitan areas.

These results suggest that urban recycling policies should not focus exclusively on expanding collection fleets, but also on strategically locating intermediate consolidation infrastructure capable of reducing collection effort in dense urban areas. The proposed framework provides quantitative support for evaluating the trade-offs between infrastructure investment, operational efficiency, and environmental impact under different urban-growth scenarios.

\section{Conclusions and Further Research}\label{sec:5}

This paper introduced the SWIFT framework, an integrated optimization approach for the design of multi-type urban waste collection and recycling systems. The proposed methodology jointly determines the locations of consolidation facilities and treatment plants, together with the associated two-echelon transportation network, while explicitly accounting for consolidation operations, multiple waste streams, and vehicle-capacity limitations. By integrating strategic infrastructure decisions with operational routing activities, the framework provides a comprehensive representation of the interactions that arise in urban recycling systems.

The computational experiments based on realistic urban scenarios derived from Medellín demonstrated the benefits of coordinated infrastructure and transportation planning. The results showed that intermediate consolidation facilities can play a key role in improving system performance by creating more compact collection regions, reducing transportation effort, and enhancing vehicle utilization. At the same time, the analysis revealed that increasing infrastructure density introduces additional complexity into the transfer network, highlighting the need to carefully balance infrastructure investments against operational gains. From a sustainability perspective, the findings indicate that integrated planning approaches can contribute to more efficient use of transportation resources and support the development of more sustainable urban recycling systems.

Several avenues for future research emerge from this study. These include the incorporation of uncertainty in waste generation and collection requirements, the consideration of dynamic and time-dependent transportation conditions, and the extension of the model to multi-period planning settings. Additional opportunities include the integration of environmental performance indicators, such as greenhouse gas emissions and energy consumption, as well as the consideration of emerging operational challenges associated with smart-city initiatives and circular-economy strategies.

{\bf Funding: } VB and YH acknowledge financial support from grants PID2020-114594GB-C21 and PID2024-156594NB-C21 (funded by MICIU/AEI/10.13039/501100011033), and from the FEDER–Junta de Andalucía project C-EXP-139-UGR23. VB also acknowledges support from the IMAG María de Maeztu grant CEX2020-001105-M. YH additionally acknowledges support from grants SOL2024-31596 and SOL2024-31708, and from the IMUS María de Maeztu grant CEX2024-001517-M (funded by MICIU/AEI/10.13039/501100011033).


\end{document}